\documentclass{amsart}

\usepackage{amssymb}
\usepackage{amsfonts, mathrsfs}
\usepackage{floatflt,graphicx}
\usepackage[all, poly, knot]{xy}
\usepackage{subfig}
\usepackage{sidecap}

\input xy
\xyoption{arc}
\xyoption{all}
\xyoption{knot}

\newtheorem{theorem}{Theorem}[section]
\newtheorem{lemma}[theorem]{Lemma}
\theoremstyle{definition}
\newtheorem{definition}[theorem]{Definition}

\theoremstyle{remark}
\newtheorem{remark}[theorem]{Remark}
\numberwithin{equation}{section}
\theoremstyle{plain}

\newtheorem{fact}[theorem]{Fact}

\newtheorem{conjecture}[theorem]{Conjecture}

\newtheorem{question}[theorem]{Question}
\newtheorem{corollary}[theorem]{Corollary}

\newtheorem{proposition}[theorem]{Proposition}

\newcommand{\cl}[1]{\text{cl}(#1)}
\newcommand{\vir}[1]{\text{vir}(#1)}
\newcommand{\tr}[1]{\text{tr}(#1)}
\newcommand{\kn}[1]{\text{kn}(#1)}
\newcommand{\un}[1]{\text{u}(#1)}

\newcommand{\viru}[1]{\text{u}_\text{v}(#1)}
\newcommand{\virtr}[1]{\text{tr}_\text{v}(#1)}
\newcommand{\ubtr}[1]{\text{tr}_\text{\"{u}}(#1)}
\newcommand{\vircl}[1]{\text{cl}_\text{v}(#1)}

\begin{document}
\title[Pseudodiagram Theory]{Classical and Virtual Pseudodiagram Theory and New Bounds on Unknotting Numbers and Genus}
\author{A. Henrich}
\address{Seattle University\\
Seattle, WA 98122}
\email{henricha@seattleu.edu}
\author{N. MacNaughton}
\address{Williams College\\
Williamstown MA, 01267}
\email{noel.f.macnaughton@williams.edu}
\author{S. Narayan}
\address{Oberlin College\\
Oberlin, OH 44074}
\email{sneha.narayan@oberlin.edu}
\author{O. Pechenik}
\address{Oberlin College\\
Oberlin, OH 44074}
\email{oliver.pechenik@oberlin.edu}
\author{J. Townsend}
\address{Scripps College \\
Claremont, CA 91711}
\email{jeniphyr@gmail.com}

\date{\today}

\begin{abstract}
A pseudodiagram is a diagram of a knot with some crossing information missing.  We review and expand the theory of pseudodiagrams introduced by R. Hanaki. We then extend this theory to the realm of virtual knots, a generalization of knots. In particular, we analyze the \emph{trivializing number} of a pseudodiagram, i.e. the minimum number of crossings that must be resolved to produce the unknot. We consider how much crossing information is needed in a virtual pseudodiagram to identify a non-trivial knot, a classical knot, or a non-classical knot. We then apply pseudodiagram theory to develop new upper bounds on unknotting number, virtual unknotting number, and genus.
\end{abstract}
\maketitle

\section{Introduction}
Recently, in \cite{hanaki}, Hanaki introduced the concept of a knot pseudodiagram and the related notions of trivializing number and knotting number.

\begin{figure}[h]
\subfloat[][]{\includegraphics[scale=.12]{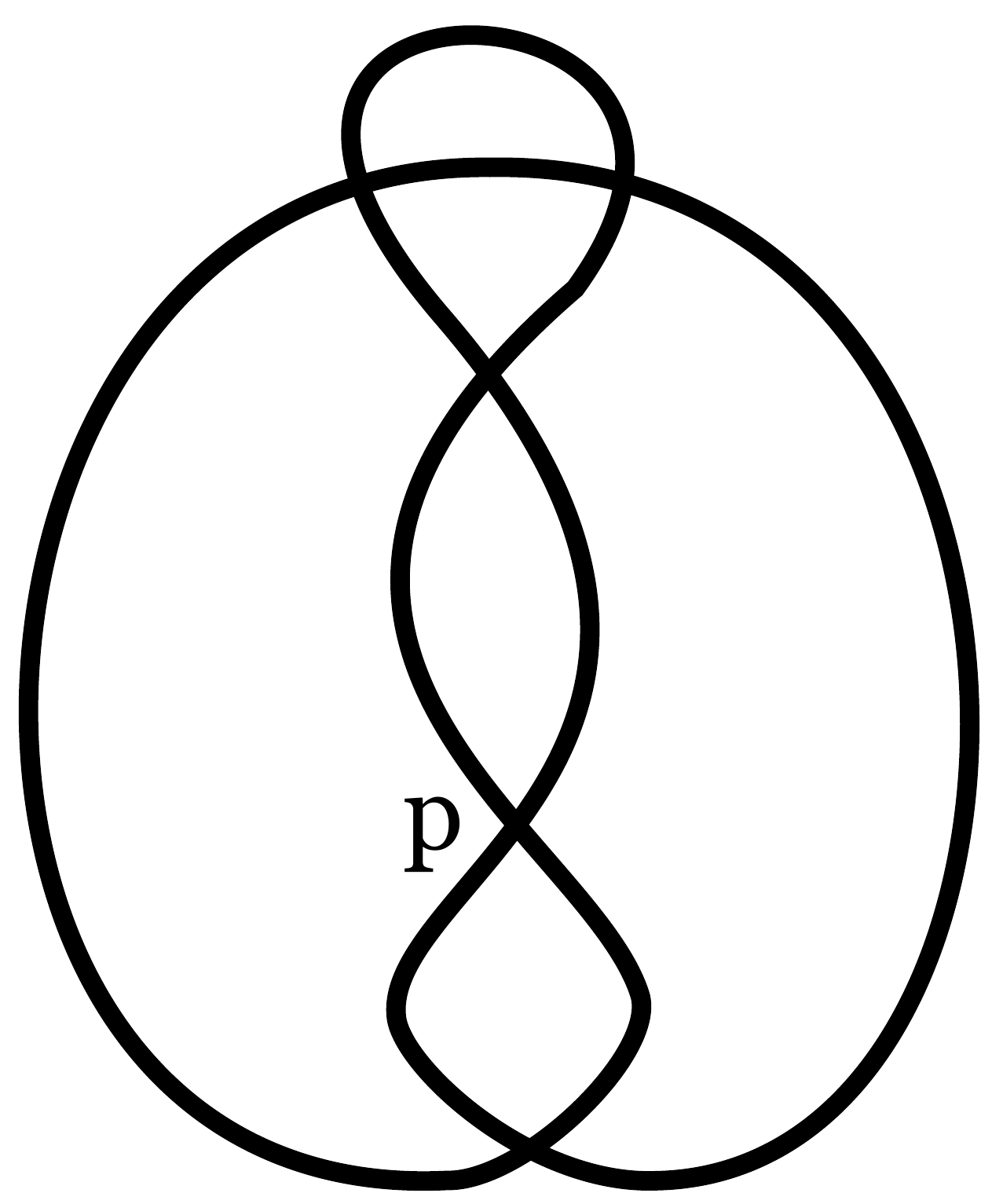}}
\subfloat[][]{\includegraphics[scale=.12]{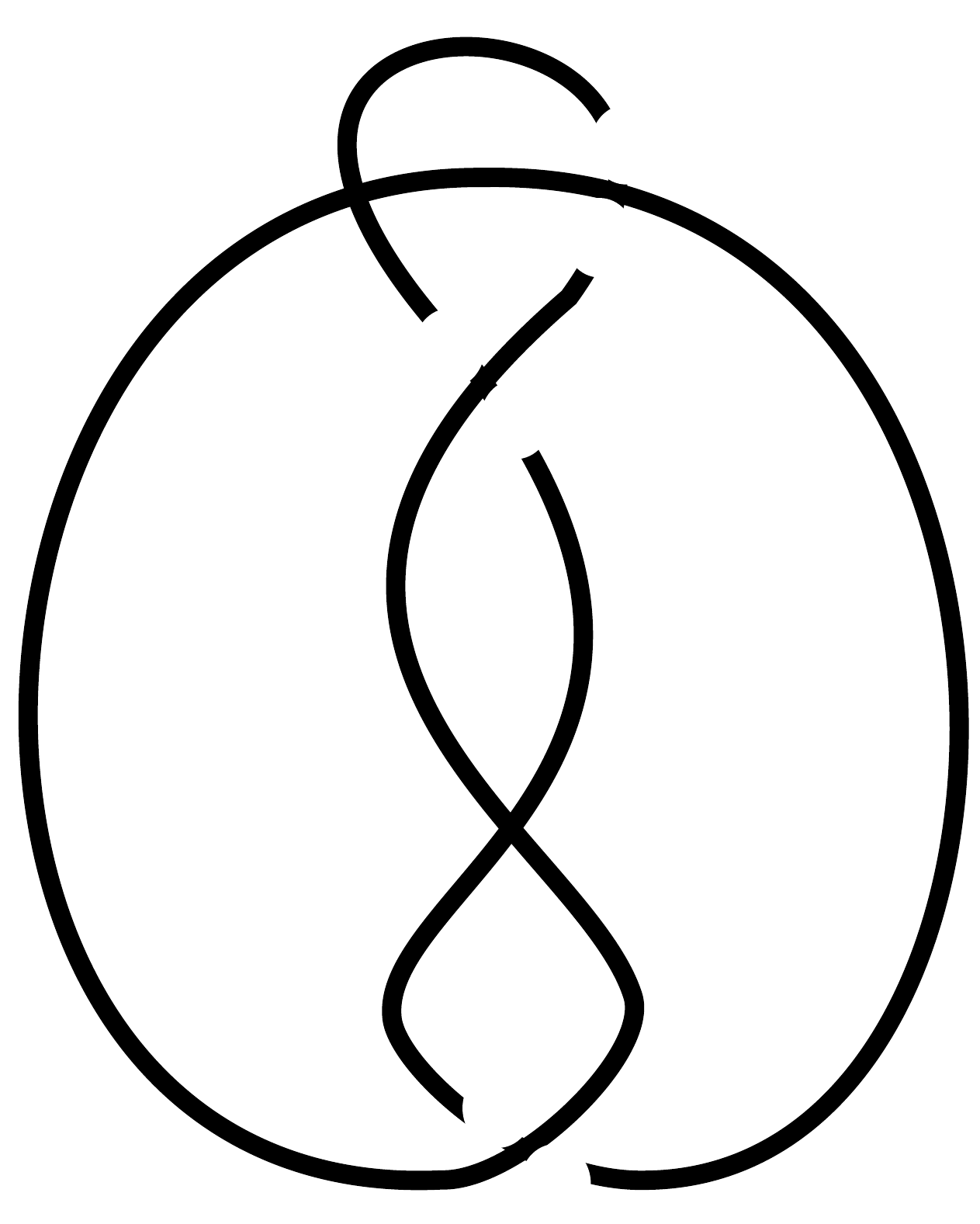}}
\subfloat[][]{\includegraphics[scale=.12]{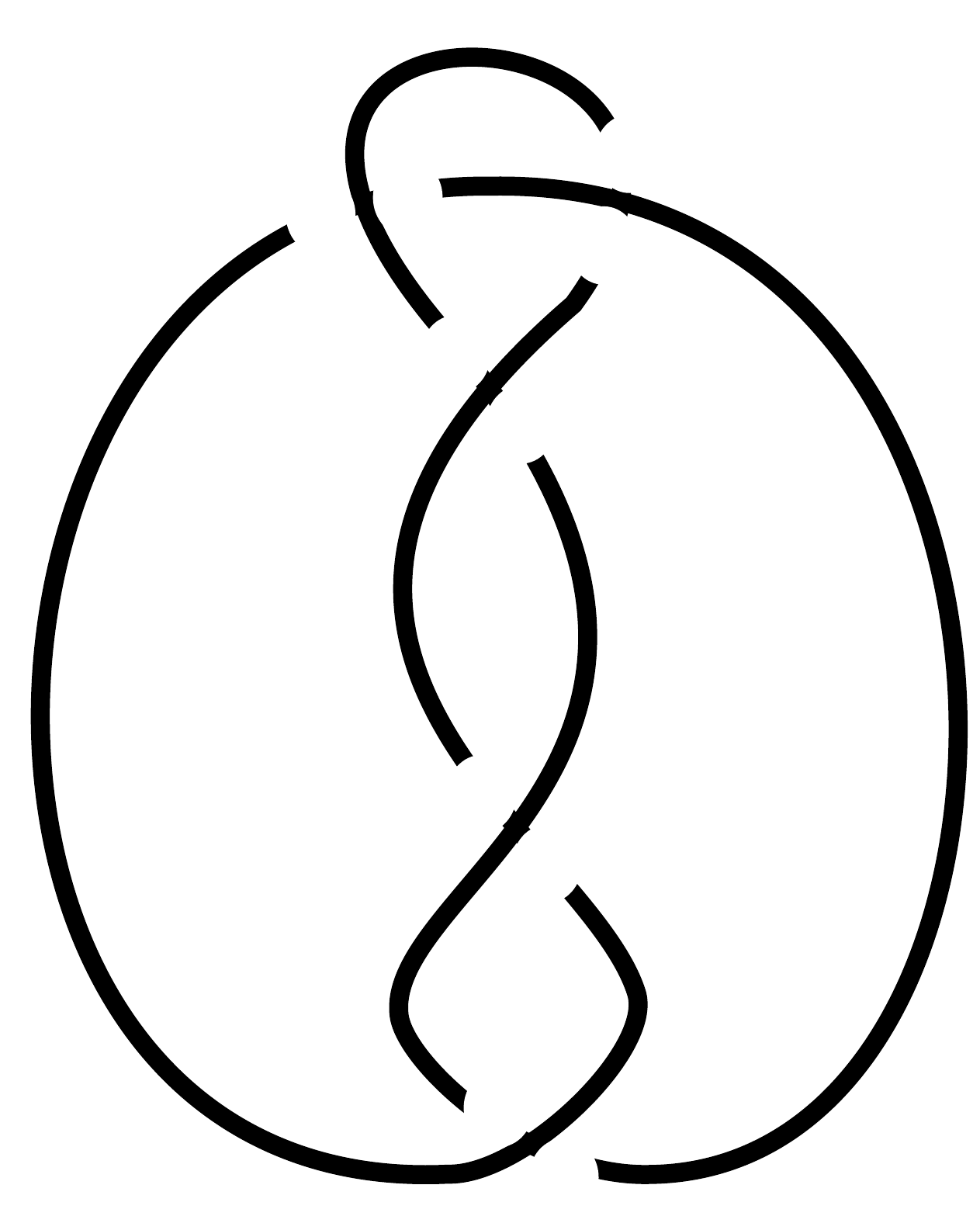}}
\caption{ \small{The three drawings above are all pseudodiagrams. Drawing (A) is a shadow where a precrossing $p$ is marked, and drawing (C) is a diagram.}}
\label{pseudofig}
\end{figure}

\begin{definition}
A \emph{pseudodiagram} $P$ is a knot diagram in which some crossings are undetermined.  Such crossings are called \emph{precrossings}.  A precrossing is represented as a flat crossing in a drawing.  We \emph{resolve} a precrossing by assigning the local writhe of that crossing.  In other words, a precrossing of a diagram is resolved by converting it to a traditional crossing.
We call a pseudodiagram  in which all crossings are undetermined a \emph{shadow} and one in which all crossings are determined a \emph{diagram}. These definitions are illustrated in Figure \ref{pseudofig}. 
\end{definition}

This concept is motivated by the study of DNA knotting.  In some pictures of DNA molecules, we are unable to determine which strand is in on top in some crossings. Pseudodiagram theory investigates what information can be determined from incomplete data about a knot's crossings.

A major concern is determining whether a pseudodiagram is necessarily knotted or unknotted regardless of how the remaining precrossings are resolved.  Hanaki introduces the following:
\begin{definition}
\label{tr}
The \emph{trivializing number} $\tr{P}$ of a pseudodiagram $P$ is the minimum number of precrossings which must be resolved so that the resulting pseudodiagram is necessarily unknotted.  If there is no resolution of precrossings such that the resulting diagram is isotopic to the unknot, then we say that $\tr{P}=\infty$. 
\end{definition}  

\begin{definition}
The \emph{knotting number} $\kn{P}$ of a pseudodiagram $P$ is the minimum number of precrossings which must be resolved so that the resulting pseudodiagram is necessarily knotted, i.e. not the unknot.  If every resolution of precrossings of $P$ results in the unknot, then we say that $\kn{P}=\infty$.
\end{definition}

In Section~\ref{Pseudodiagrams}, we discuss values of these quantities. We provide prerequisite background from classical knot theory in Section~\ref{Prelims} to aid this discussion.
 
In Section~\ref{Virtual Pseudodiagrams}, we define a generalization of knot theory introduced by Kauffman in~\cite{VKT} called virtual knot theory. We extend the notion of a pseudodiagram to this broader class of knots in two natural ways. In Sections~\ref{Results on Characteristic Values of Pseudodiagrams} and~\ref{Kauffman's $J$-invariant}, we explore results pertaining to these virtual pseudodiagrams.
 
Finally, in Section~\ref{Unknotting Numbers and Genus}, we use pseudodiagram theory to find bounds for the classical and virtual unknotting numbers as well as the canonical genus for classical knots.

\section{Preliminaries}
\label{Prelims}
We define a \emph{knot} to be an embedded copy of the circle in $\mathbb{R}^3$ up to ambient isotopy.  A \emph{link} is an embedded copy of one or more circles in $\mathbb{R}^3$. Note that we may alternatively consider knots and links as sitting inside a thickened 2-sphere. When we study knots, we typically consider diagrams of knots---projections of knots onto a two-dimensional plane. These diagrams are pictured generically as curves with no self-tangencies or triple-points, only crossings given as double-points decorated to show which strand of the knot passes over and which passes under at the crossing. Two knot diagrams represent equivalent knots if and only if the diagrams can be related by a sequence of Reidemeister moves and ambient isotopy in the plane. 

\begin{figure}[h]
\centering
\includegraphics[width=.2\textwidth]{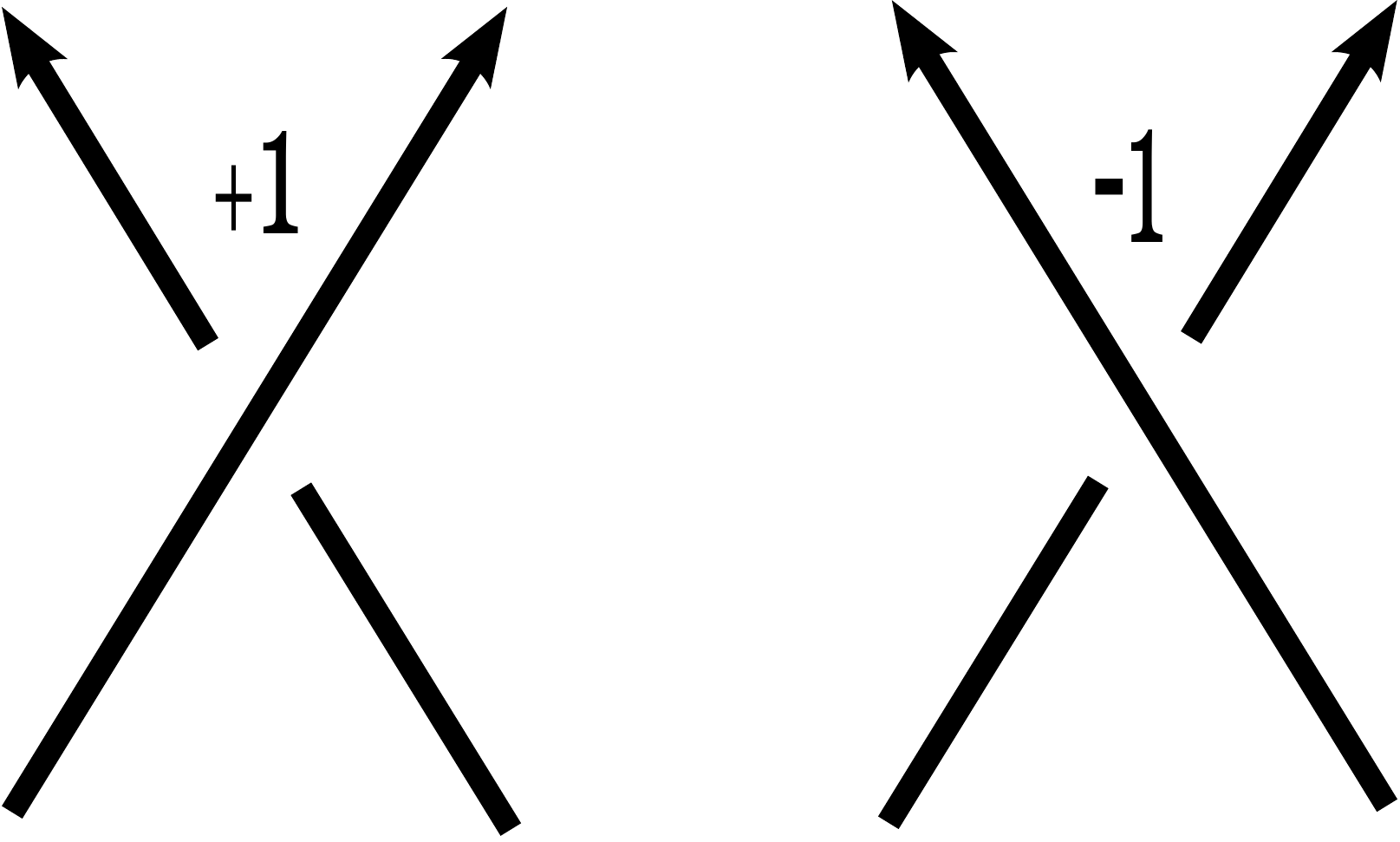}
\caption{\small{Local writhes of crossings.}}\label{Signs}
\end{figure}

Another way of representing oriented knots is using objects called \emph{Gauss diagrams}, which can be associated in a natural way to knot diagrams. A Gauss diagram is a counter-clockwise oriented circle (the `core' circle) parametrizing the knot equipped with signed oriented chords. The chords represent crossings in the knot diagram, while the signs and orientations contain information about under- and over-strands. To be more precise, a chord is oriented from the preimage of the over-strand of the associated crossing to the preimage of the under-strand on the core circle. Each chord is given the sign of the local writhe of its corresponding crossing in the knot diagram, as illustrated in Figure~\ref{Signs}. Figure~\ref{Gauss} gives an example of a knot diagram and its associated Gauss diagram. 

A \emph{base point} for a Gauss diagram is a point on the core circle at which we begin to traverse the core circle with respect to the orientation.  This base point corresponds to a base point on the knot.

\def\TrefoilA{\xygraph{!{0;/r2.0pc/:}
!P3"a"{~>{}}
!P9"b"{~:{(1.3288,0):}~>{}}
!P3"c"{~:{(2.5,0):}~>{}}
!{\vover~{"b2"}{"b1"}{"a1"}{"a3"}=>|{1}}
!{"b4";"b2" **\crv{"c1"}}
!{\vover~{"b5"}{"b4"}{"a2"}{"a1"}=>|{3}}
!{"b7";"b5" **\crv{"c2"}}
!{\vover~{"b8"}{"b7"}{"a3"}{"a2"}=>|{2}}
!{"b1";"b8" **\crv{"c3"}}}}

\def\TrefoilGauss{
\begin{xy} /r15mm/:
,{\ellipse<>{}}
,(2,0)="1" ,*+!L{1}
,(.292893, .707107)="3" ,*+!DR{3}
,(1,-1)="2" ,*+!U{2}
,(0,0)="a1" ,*+!LU{+}
,(1.707107, -.707107)="a3" ,*+!D{+}
,(1,1)="a2" ,*+!LU{+}
,{\ar@{->} 0; "1"}
,{\ar@{->} (1,1); (1, -1)}
,{\ar@{<-} (.292893, .707107); (1.707107, -.707107)}
\end{xy}
}

\begin{figure}[h]
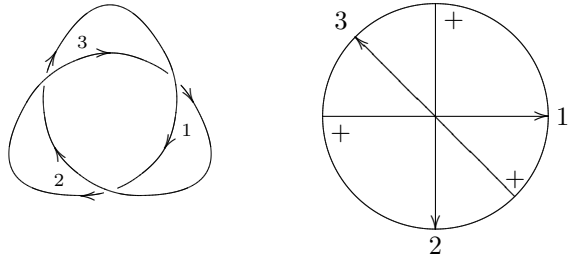

$$\TrefoilA\hspace{10mm}\TrefoilGauss$$
\caption{\small{A trefoil knot and its Gauss diagram.}}\label{Gauss}
\end{figure}

There is a notion of equivalence for Gauss diagrams that corresponds to the equivalence of knot diagrams. In particular, there are combinatorial Reidemeister-type moves that relate Gauss diagrams representing the same knot, as described in~\cite{PolyakViro}. 

We will especially make use of two classes of figures related to Gauss diagrams. An \emph{arrow diagram} is a Gauss diagram that is missing sign information on its chords. Furthermore, a \emph{chord diagram} is an unoriented arrow diagram. Hence, a chord diagram preserves only the core circle and the chords from the original Gauss diagram. Chord diagrams may also be thought of as free knots. See~\cite{freeknot} and~\cite{freelink} for more on free knots. From an arrow diagram or a chord diagram, it is not possible to uniquely reconstruct the original knot diagram.  More information about Gauss diagrams, arrow diagrams, and chord diagrams can be found in~\cite{manturov}.

\begin{definition}
We call two chords of a Gauss, arrow, or chord diagram \emph{parallel} if they can be drawn so that they do not intersect.  We say that a Gauss, arrow or chord diagram is \emph{parallel} if all the chords are pairwise parallel.
\end{definition}

\section{Results for Pseudodiagrams}
\label{Pseudodiagrams}
In the introduction, we discussed the notion of a pseudodiagram as well as the notions of trivializing and knotting number. Here, we present what is known about values of these numbers.

\begin{theorem}[Hanaki, \cite{hanaki}] 
The trivializing number of any shadow is even. 
\label{hanaki}
\end{theorem}

In the section below, we modify Hanaki's proof and prove Theorem \ref{H2}, a generalization of Theorem \ref{hanaki}.

\begin{definition}
A set $\mathscr{T}$ of precrossings of a pseudodiagram $P$ that can be resolved so that $P$ is necessarily the unknot is called a \emph{trivializing set} of $P$.  

We can further define two types of trivializing sets:
\begin{description}
\item[basic] If there exists no proper subset $\mathscr{U} \subsetneq \mathscr{T}$ such that $\mathscr{U}$ is a trivializing set of $P$, then $\mathscr{T}$ is a \emph{basic trivializing set} of $P$.  

\item[minimum] If $\mathscr{T}$ is a trivializing set such that $|\mathscr{T}|\le |\mathscr{U}|$ for all trivializing sets $\mathscr{U}$, then $\mathscr{T}$ is a \emph{minimum trivializing set} of $P$.
\end{description}
\end{definition}

Observe that $\tr{P}=|\mathscr{T}|$ where $\mathscr{T}$ is a minimum trivializing set, and also that not every basic trivializing set of crossings is minimum.
\begin{theorem} \label{H2}
The cardinality of every basic trivializing set of a shadow is even. 
\end{theorem}

In order to prove this theorem, we require the following additional lemmas.


\begin{lemma}\label{||triv}
A chord diagram with only parallel chords is trivial.
\end{lemma}
\begin{proof}
Triviality is shown by applying Reidemeister move I (for chord diagrams), once for each chord of the chord diagram.
\end{proof}
\begin{lemma}\label{+nontriv}
If the (classical) pseudodiagram $P$ contains exactly two precrossings and the corresponding chords intersect in the chord diagram of $P$, then $P$ can be resolved nontrivially.
\end{lemma}
\begin{proof}
We first interpret the second degree Vassiliev invariant $\nu_2$ in terms of Gauss diagrams with our precrossings $a$ and $b$ as singular crossings.  For an introduction to Vassiliev invariants and singular knots, see \cite{Goussarovpaper} and  \cite{Vassiliev}.

As was illustrated by Polyak and Viro in \cite{PolyakViroGD}, and Chmutov, Khoury and Rossi in \cite{chmutov}, the second degree Vassiliev invariant $\nu_2$ is defined on realizable Gauss diagrams as  $$\nu_2(K)=\sum_{(x,y)\in C} w_xw_y,$$ where $w_x$ is the local writhe at crossing $x$. Given an arbitrary base point on the core circle of the Gauss diagram of $P$, we include an ordered pair of chords $(x,y)$ in $C$ if, on one counterclockwise circuit around the core circle of the Gauss diagram of $P$, we encounter $x$ and $y$ chord endpoints in the order $x_H,y_T,x_T,y_H$, where $x_H$ denotes the head of the $x$ chord arrow, and $y_T$ denotes the tail of the $y$ chord arrow.  Possible relations of chord pairs $(x,y)$ are illustrated in Figure \ref{ijij}. Only the chord pair labeled (a) in Figure \ref{ijij} is in $C$.  Although the set $C$ depends on the base point, $\nu_2$ is independent of the choice of base point. 

It should be noted that $\nu_2$ evaluated on the unknot is 0.

\begin{figure}[h]
\includegraphics[scale=.35]{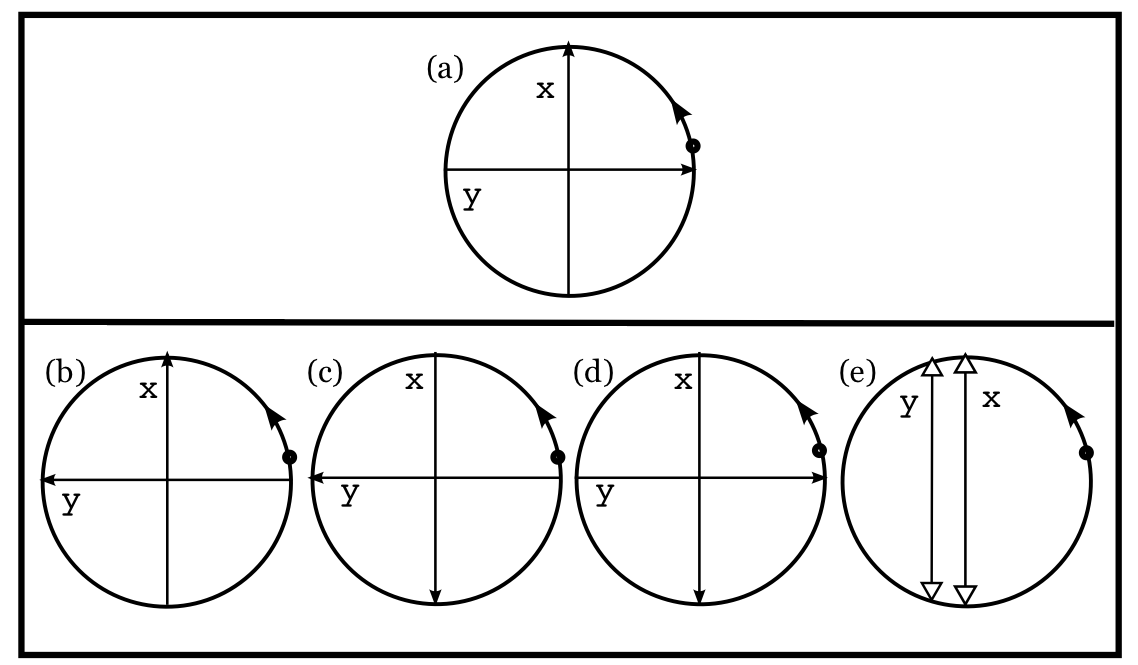}
\caption{\small{The crossing pair $(x,y)$ can be related in the chord diagram as shown in the figure above: there are four crossing relations, and four parallel relations. Only the chord relation (a) is in $C$, and therefore is counted for $\nu_2$.}}
\label{ijij}
\end{figure}
Consider the two precrossings $a$ and $b$ as singular crossings.  Let the base point for defining $C$ directly follow an endpoint of $b$ (as in Figure \ref{CD_nontriv}), and consider all chords relative to the quadrants pictured in Figure \ref{CD_nontriv}. We define $P_{HT}$ to be the diagram obtained from $P$ by resolving $a$ so that we meet the head of $a$ first when traversing the core circle, and by resolving $b$ so that we meet its tail first. We define $P_{HH}, P_{TH},$ and $P_{TT}$ similarly.  Then, applying the second derivative of Vassiliev invariants on singular knots, $$\nu_2(P)=\nu_2(P_{HH})+\nu_2(P_{TT})-\nu_2(P_{HT})-\nu_2(P_{TT}).$$

We wish to show that there are an odd number of chord pairs $(x,y)$ which contribute to $\nu_2(P)$. This will guarantee that $\nu_2(P)\ne 0$, as every chord pair in $C$ contributes either $+1$ or $-1$.  Since $0\ne \nu_2(P)=\nu_2(P_{HH})+\nu_2(P_{TT})-\nu_2(P_{HT})-\nu_2(P_{TT})$, it follows that there exists a resolution of the precrossings $a$ and $b$ such that the resulting knot ($P_{HH}, P_{TT},P_{TH},$ or $P_{HT}$) has $\nu_2\ne 0$, and therefore is nontrivial.

\begin{SCfigure}
\centering
\includegraphics[trim= 100 500 100 0, width=35mm]{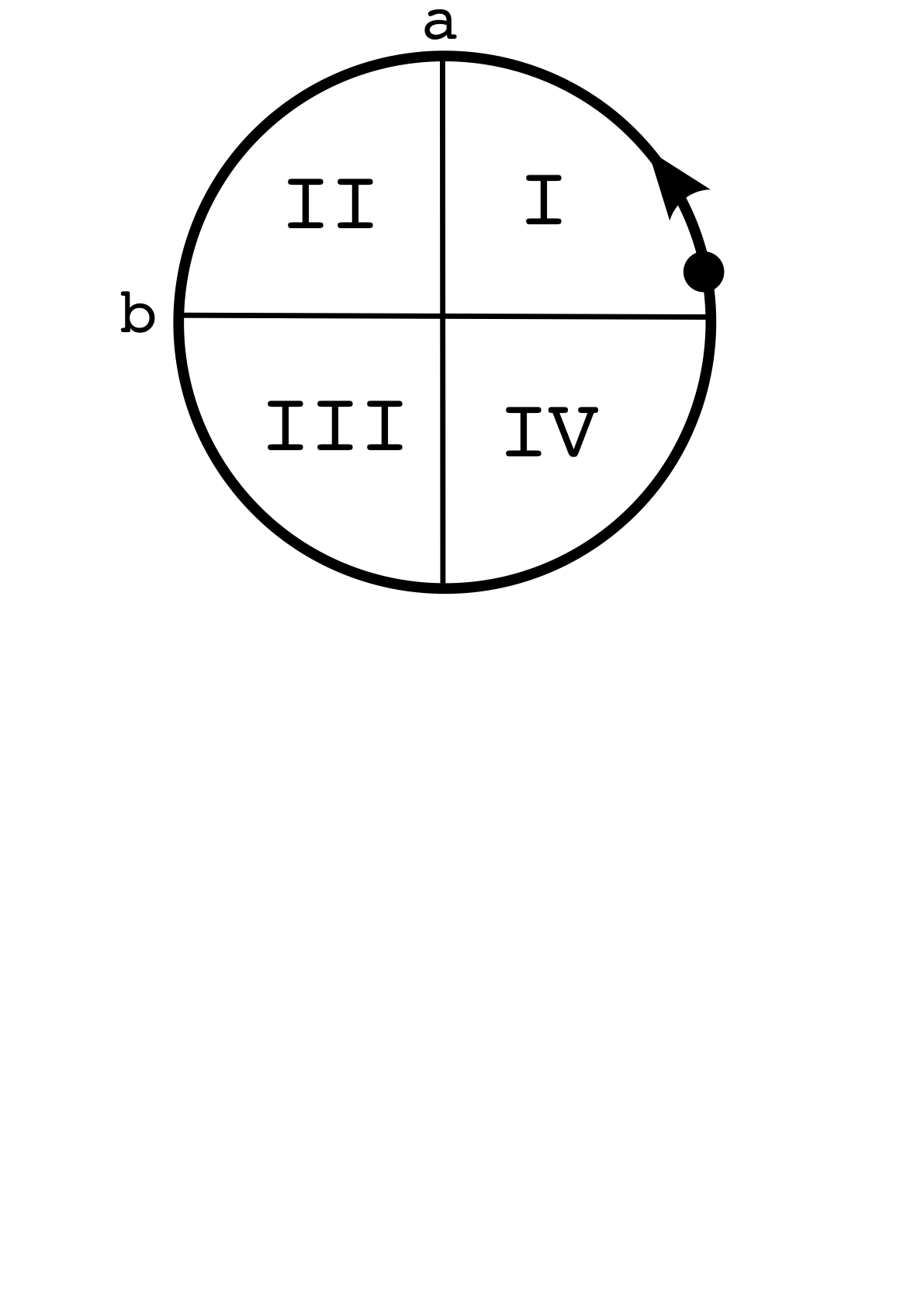}

\caption{\small{The chords $a$ and $b$ correspond to precrossings in $P$, and divide the chord diagram into four quadrants.}}

\label{CD_nontriv}
\end{SCfigure}

Observe that the resolution of precrossings $a$ and $b$ does not affect whether a chord pair $(x,y)$ (where $x,y\notin \{a,b\}$) is in $C$.  Such a chord pair $(x,y)$ is either counted for computations of \emph{all} or \emph{none} of $\nu_2(P_{HH}), \nu_2(P_{HT}), \nu_2(P_{TT}),\nu_2(P_{TH})$. Therefore $(x,y)$ appears either 4 or 0 (both even) times in the calculation of $\nu_2(P)$. Thus to calculate $\nu_2(P)$ over $\mathbb{Z}_2$ it suffices to consider only chord pairs $(x,y)\in C$ where either $x$ or $y$ is $a$ or $b$. 

A chord $x$ intersects $a$ or $b$ if and only if the endpoints of $x$ are in different quadrants of Figure \ref{CD_nontriv}.  The statement that $(x,b)\in C$ depends on how $b$ is oriented and which quadrants the heads and tails of $c$ are in.  We can classify the number of times any chord from quadrant $A$ to quadrant $B$ appears in the calculation of $\nu_2(P_{HH}),\nu_2(P_{TT}),\nu_2(P_{HT}),$ and $\nu_2(P_{TH})$. We represent this information below as the matrices $HH$, $TT$, $HT$ and $TH$, where $HH_{ij}$ gives the number of times any chord $x_{i\to j}$ with tail in quadrant $i$ and head in quadrant $j$ is counted in a chord pair with $a$ or $b$ in $C$ for the computation of $\nu_2(P_{HH})$.
$$
HH=
\begin{pmatrix}
0 & 0 & 0 & 0\\
0 & 0 & 0 & 1\\
0 & 0 & 0 & 1\\
0 & 0 & 1 & 0\\
\end{pmatrix}
\hspace{1in}
TT=\begin{pmatrix}
0 & 0 & 0 & 0 \\
1 & 0 & 0 & 0\\
2 & 1 & 0 & 0\\
1 & 1 & 0 & 0\\
\end{pmatrix}$$ 

$$
HT=\begin{pmatrix}
0 & 0 & 0& 0\\
0 & 0 & 0 & 1\\
1 & 1 & 0 & 1\\
1 & 1 & 1 & 0\\
\end{pmatrix}
\hspace{1in}
TH=\begin{pmatrix}
0 & 0 & 0& 0\\
1 & 0 & 0 & 0\\
1 & 0 & 0 & 0\\
0 & 0 & 0 & 0\\
\end{pmatrix}
$$
 Given these matrices, the $ij$ entry of the sum matrix $HH+TT+HT+TH$ gives the total number of times any chord $x_{i\to j}$ appears in the calculation of $\nu_2(P)$. 
 $$
HH+TT+HT+TH=
\begin{pmatrix}
0 & 0 & 0 & 0\\
2 & 0 & 0 & 2\\
4 & 2 & 0 & 2\\
2 & 2 & 2 & 0\\
\end{pmatrix}
 $$
 From this it is clear that every chord pair $(x,a),(a,x),(x,b),$ or $(b,x)$ appears an even number of times when we calculate $\nu_2(P)$.  We have considered every chord pair except for $(a,b)$ and $(b,a)$, therefore the parity of the number of $w_xw_y$ terms summed is precisely the parity of chord pair terms $(a,b)$ and $(b,a)$ in our four resolutions of $a$ and $b$.  Observe that $(b,a) \notin C$ for any resolutions of $b$ and $a$, and $(a,b)$ satisfies the first chord crossing relation if and only if $a$ and $b$ are resolved as in $P_{HT}$.  Thus the total number of contributing chord pair terms is odd, and $\nu_2(P)\ne 0$.
\end{proof}

\begin{lemma}\label{+nontriv2}
Any pseudodiagram $P$ with a chord diagram containing intersecting prechords (i.e. prechords correspond to precrossings of $P$) can be resolved nontrivially.
\end{lemma}
\begin{proof}
Arbitrarily resolve all but two precrossing $a$ and $b$ of $P$ such that the chords of $a$ and $b$ intersect in a chord diagram.  By Lemma \ref{+nontriv}, the resulting pseudodiagram has a nontrivial resolution.  Thus $P$ has a nontrivial resolution.
\end{proof}

\begin{lemma}\label{triv=crosses}
The trivializing number of a shadow $S$ is equal to the minimum number of chords that can be deleted from the chord diagram such that all remaining chords are parallel.
\end{lemma}
\begin{proof}
Let $\mathscr{T}$ be a trivializing set of precrossings for $S$.  By Lemma \ref{+nontriv2}, since $\mathscr{T}$ is a trivializing set, deleting the corresponding chords in the chord diagram of $S$ must leave only parallel chords.  Therefore, the trivializing number is at least the minimum number of chords that can be deleted to leave only parallel chords.

We will call a chord $c$ an \emph{exterior chord} if it cuts off an arc with only endpoints of chords in $\mathscr{T}$.  Consider the strand $a$ of $S$ corresponding to this arc. Any precrossings along $a$ must be in $\mathscr{T}$.  By resolving the precrossings on $a$ so that $a$ lies over all other strands of $S$, we may contract $a$ so as to remove all crossings on it. This allows us to perform a Reidemeister I move that eliminates the crossing $c$. 

This produces a new shadow with fewer crossings that can be trivialized by the remaining elements of $\mathscr{T}$. Hence, we may iterate this process on an exterior chord of the resulting shadow until all crossings and precrossings of $S$ are eliminated, resulting in a diagram of the unknot. This process can be performed with any trivializing set---in particular with a minimum trivializing set. Thus, the trivializing number is equal to the minimum number of chords that can be deleted to leave only parallel chords.
\end{proof}

\begin{proof}[Proof of Theorem \ref{H2}.]
Let $S$ be a shadow, and $\mathscr{T}$ a basic trivializing set of precrossings for $S$. We wish to show that $|\mathscr{T}| $ is even.  This theorem is a corollary of Lemma \ref{triv=crosses}.

It is well-known that every knot diagram is evenly intersticed (see for example \cite{VKT}).\footnote{This fact does not extend to virtual knot diagrams, a concept we will discuss in Section \ref{Virtual Pseudodiagrams}.}  In terms of chord diagrams, this means the arc cut off by every exterior chord, as described in the proof of Lemma \ref{triv=crosses}, contains an even number of chord endpoints.  All of these endpoints are from chords in $\mathscr{T}$.  Because $\mathscr{T}$ is basic, each chord in $\mathscr{T}$ contributes at most one of these endpoints. 

Thus, in the proof of Lemma \ref{triv=crosses}, to contract $a$, we resolved and then eliminated an even number of crossings in $\mathscr{T}$.  Therefore, $|\mathscr{T}|$ is even.
\end{proof}

\begin{remark} We note here that it may be possible to strengthen these results regarding trivializing numbers using Manturov's work on free knots and links. We hope this will be a topic of future study.
\end{remark}

In \cite{hanaki}, Hanaki proved that $\kn{S}\ge 3$ for any shadow $S$.  We offer an alternative proof of this fact.
\begin{proposition}\label{knottingnotsmall}
For any shadow $S$, $\kn{S}\ge 3$.
\end{proposition}
\begin{proof}
 Suppose that we resolve only the precrossing $c_1$. We traverse $S$ from a base point $b$ such that we first encounter $c_1$ as an overpass, resolving all precrossings as overpasses at the first encounter.    It is well known (see, for example, \cite{knotbook}) that this resolution is the unknot.
 
Now suppose we resolve the precrossings $c_1$ and $c_2$.  The base point $b$ is arbitrary, and there exists a base point $\hat{b}$ and an oriented traversal such that we first encounter both $c_1$ and $c_2$ as overpasses.  Again, resolving all precrossings along this traversal as overpasses at the first encounter produces the unknot.
\end{proof}



\section{Virtual Pseudodiagrams}
\label{Virtual Pseudodiagrams}
In this section, we extend Hanaki's notion of a pseudodiagram into the domain of virtual knots, a generalization of knots.  Virtual knots can be interpreted as knots on surfaces of various genera, so they serve as a model to examine the knottedness of biological polymers wrapped around cellular structures.  Before introducing virtual pseudodiagrams, we review certain pertinent definitions from virtual knot theory. Additional introductory material on virtual knots may be found in \cite{VKT}.

\begin{figure}
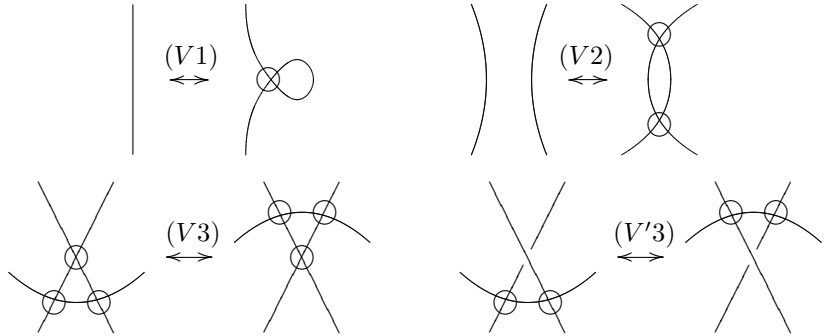


\[
\xy
(-30,-10)*{}; (-30,10)*{} **\dir{-};
{\ar@{<->}(-25,0)*{}; (-20,0)*{}}; ?(.75)*dir{>}+(-1, 3)*{(V1)};
(-15, 10)*{}; (-15, -10)*{} **\crv{(-15,5) & (-12,0) & (-8,-4) & (-5,0) & (-8, 4) & (-12, 0) & (-15, -5)};
(-12,0)*{\bigcirc};
(15,10)*{}="E";
(15, -10)*{}="F";
"E"; "F" **\crv{(19, 0)};
(25,10)*{}="G";
(25, -10)*{}="H";
"G"; "H" **\crv{(21, 0)};
{\ar@{<->}(33,0)*{}; (28,0)*{}}; ?(.75)*dir{>}+(1, 3)*{(V2)};
(35,10)*{}="I";
(35,-10)*{}="J";
(45, 10)*{}="K";
(45, -10)*{}="L";
"I"; "J" **\crv{(40,7) & (43,0) & (40, -7)};
"K"; "L" **\crv{(40, 7) & (37, 0) & (40, -7)};
(40,6)*{\bigcirc};
(40,-6)*{\bigcirc}; 
\endxy
\]

\[
\xy
(75,10)*{}="AT";
(85, -10)*{}="AB";
(85, 10)*{}="BT";
(75, -10)*{}="BB";
(71, -2)*{}="CL";
(89, -2)*{}="CR";
"CL"; "CR" **\crv{(80, -10)};
"AT"; "AB" **\crv{};
"BT"; "BB" **\crv{};
(80,0)*{\bigcirc};
(77,-6)*{\bigcirc};
(83,-6)*{\bigcirc};
{\ar@{<->}(92,0)*{}; (98,0)*{}}; ?(.75)*dir{>}+(-1, 3)*{(V3)};
(105,10)*{}="A'T";
(115, -10)*{}="A'B";
(115, 10)*{}="B'T";
(105, -10)*{}="B'B";
(101, 2)*{}="C'L";
(119, 2)*{}="C'R";
"C'L"; "C'R" **\crv{(110, 10)};
"A'T"; "A'B" **\crv{};
"B'T"; "B'B" **\crv{};
(110,0)*{\bigcirc};
(107,6)*{\bigcirc};
(113,6)*{\bigcirc};
(135,10)*{}="DT";
(145, -10)*{}="DB";
(145, 10)*{}="ET";
(135, -10)*{}="EB";
(131, -2)*{}="FL";
(149, -2)*{}="FR";
"FL"; "FR" **\crv{(140, -10)};
"DT"; "DB" **\crv{}; \POS?(.5)*{\hole}="a";
"ET"; "a" **\crv{};
"a"; "EB" **\crv{};
(137,-6)*{\bigcirc};
(143,-6)*{\bigcirc};
{\ar@{<->}(152,0)*{}; (158,0)*{}}; ?(.75)*dir{>}+(-1, 3)*{(V'3)};
(165,10)*{}="D'T";
(175, -10)*{}="D'B";
(175, 10)*{}="E'T";
(165, -10)*{}="E'B";
(161, 2)*{}="F'L";
(179, 2)*{}="F'R";
"F'L"; "F'R" **\crv{(170, 10)};
"D'T"; "D'B" **\crv{}; \POS?(.5)*{\hole}="a'";
"E'T"; "a'" **\crv{};
"a'"; "E'B" **\crv{};
(167,6)*{\bigcirc};
(173,6)*{\bigcirc};
\endxy
\]
\caption{\small{The Reidemeister moves for virtual knots. Two virtual diagrams are equivalent if and only if one can be obtained from the other by a sequence of these moves and the classical Reidemeister moves.}}
\label{rm}
\end{figure}

\begin{definition}
A \emph{virtual knot} is an equivalence class of knot diagrams with an additional crossing type. Instead of requiring all crossings to be either positive or negative, we introduce a third possibility, which we call a \emph{virtual crossing} and denote by drawing a small circle around the crossing. Two such knot diagrams are considered equivalent if one can be transformed into the other by a sequence of classical and/or virtual Reidemeister moves (shown in Figure \ref{rm}).  
\end{definition}

The virtual Reidemeister moves are equivalent to the \emph{virtual detour move}, which states that any strand with fixed endpoints and no classical crossings can be replaced with any other strand with the same endpoints and no classical crossings.

\begin{definition}
The minimum number of virtual crossings over all diagrams of the knot $K$ is the \emph{virtual crossing number} $c_v(K)$. 
\end{definition} 
If $c_v(K)=0$, we call $K$ a \emph{classical knot}.  If $c_v(K)\ge 1$, we say $K$ is non-classical. See~\cite{Arrow} or~\cite{am} for more on the virtual crossing number.

\begin{figure}
\includegraphics[width=.8\textwidth]{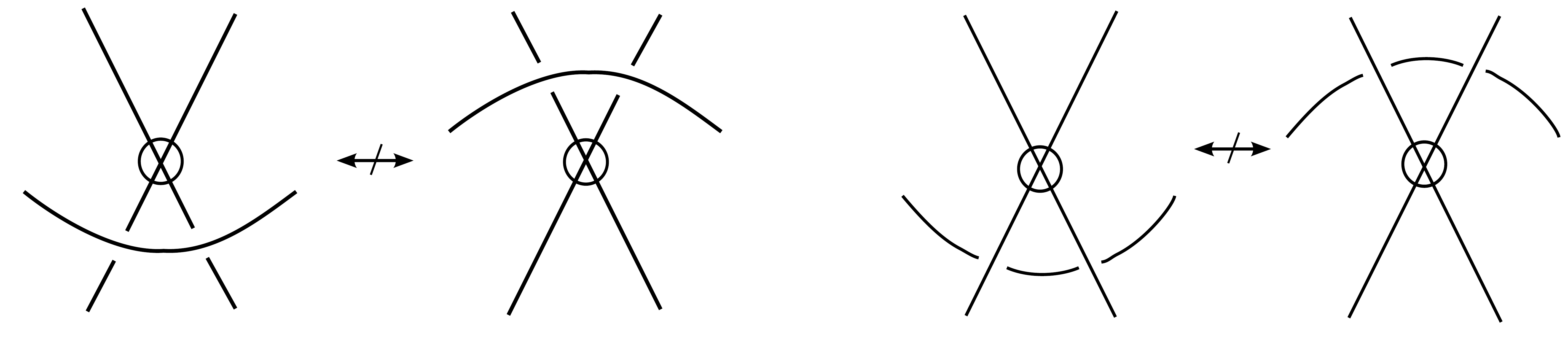}
\caption{\small{These Reidemeister-like moves are {\bf forbidden}. Allowing these two moves would render all virtual knots equivalent.} }
\label{forb}
\end{figure}

It is important to note that the two Reidemeister-like moves shown in Figure \ref{forb} are forbidden.  As shown in \cite{forbidden2} and \cite{forbidden}, allowing one of these forbidden moves leads to the theory of \emph{welded knots}, and allowing both forbidden moves trivializes the theory: all knots become equivalent to the unknot.

The motivations behind virtual knot theory are more fully explained in \cite{VKT}. Here, we merely note two equivalent definitions. Recall that classical knots may be represented as Gauss diagrams. However, not every Gauss diagram is realizable as a diagram of a classical knot.  On the other hand, every Gauss diagram has a realization as a \emph{virtual} knot.  (Virtual crossings are not recorded on the Gauss diagram.) Therefore, we may equivalently regard a virtual knot as an equivalence class of Gauss diagrams under appropriate Gauss diagram formulations of the Reidemeister moves.

Another formulation of virtual knot theory, proven to be equivalent in  \cite{CKS}, \cite{Kuperberg}, is to regard virtual knots as equivalence classes of embeddings of a circle into a thickened surface $M$. In this formulation, classical knots are precisely the equivalence classes of embeddings of circles into the thickened sphere. Non-classical virtual knots are the equivalence classes of embeddings of circles on thickened surfaces of strictly higher topological genus. Thinking geometrically, virtual crossings are artifacts of projections. Strands may not cross on $M$, but may cross in a projection of $M$ onto a plane. For every virtual knot $K$, there is a minimum topological genus of surfaces on which $K$ can be drawn.

There are two meaningful ways to extend the concept of pseudodiagrams to virtual knots.
\begin{definition}
A \emph{virtual pseudodiagram} is a diagram of a virtual knot in which some classical crossings are undetermined (however, all virtual crossings are given).  These precrossings can be resolved as positive or negative \textbf{classical} crossings, but not as virtual crossings. A \emph{virtual shadow} is a virtual pseudodiagram in which all classical crossings are undetermined.
\end{definition}

\begin{definition}
An \emph{\"{u}ber-virtual pseudodiagram} is a diagram of a virtual knot in which some crossings are undetermined precrossings. (Crossings in such a diagram may be virtual, positive classical, negative classical or precrossings.) These precrossings can be resolved as \textbf{virtual}, or positive/negative classical crossings.
\end{definition}

\begin{definition}
A \emph{virtual shadow} is either a virtual pseudodiagram or an \"{u}ber-virtual pseudodiagram, in which no classical crossing information is given.
\end{definition}

In classical pseudodiagrams, we considered the trivializing and knotting numbers.  These notions extend naturally to virtual and \"{u}ber-virtual pseudodiagrams.  In addition, two further numbers are significant:
\begin{definition}
The \emph{classicalizing number} $\cl{V}$ of a virtual pseudodiagram $V$ is the minimum number of undetermined crossings which must be (classically) resolved such that the resulting virtual pseudodiagram is necessarily classical regardless of how the remaining undetermined crossings are resolved.  If there is no resolution of undetermined crossings such that the resulting diagram is isotopic to a classical knot, then we say that $\cl{V}=\infty$.
\end{definition}

\begin{definition}
The \emph{virtualizing number} $\vir{V}$ of a virtual pseudodiagram $V$ is the minimum number of undetermined crossings which must be (classically) resolved such that the resulting pseudodiagram is necessarily non-classical regardless of how the remaining undetermined crossings are resolved.  If there is no resolution of undetermined crossings such that the resulting diagram is not isotopic to any classical knot, then we say that $\vir{V}=\infty$.
\end{definition} 

We also consider variations on the classical trivializing number.
 
\begin{definition} 
\label{virtr}
The \emph{virtual trivializing number} $\virtr{\ddot{U}}$ of an \"{u}ber-virtual pseudodiagram $\ddot{U}$ is the minimum number of precrossings which need to be determined as virtual so that the resulting \"{u}ber-virtual pseudodiagram is necessarily unknotted.
\end{definition}

\begin{definition} 
\label{ubertr}
The \emph{\"{u}ber trivializing number} $\ubtr{P}$ of an \"{u}ber-virtual pseudodiagram $P$ is the minimum number of precrossings which must be determined as either classical (positive or negative) or virtual so that the resulting \"{u}ber-virtual pseudodiagram is necessarily unknotted.\
\end{definition}

\section{Results on Characteristic Values of Pseudodiagrams}
\label{Results on Characteristic Values of Pseudodiagrams}

As a trivial knot is always classical and a non-classical knot is always knotted, the following fact is immediate. 

\begin{fact}
\label{ingies} 
For any virtual pseudodiagram $P$, $$\cl{P} \leq \tr{P} \text{   and    } \kn{P} \leq \vir{P}.$$ 
\end{fact}

\begin{figure}[h]
\subfloat[][]{\includegraphics[scale=.22]{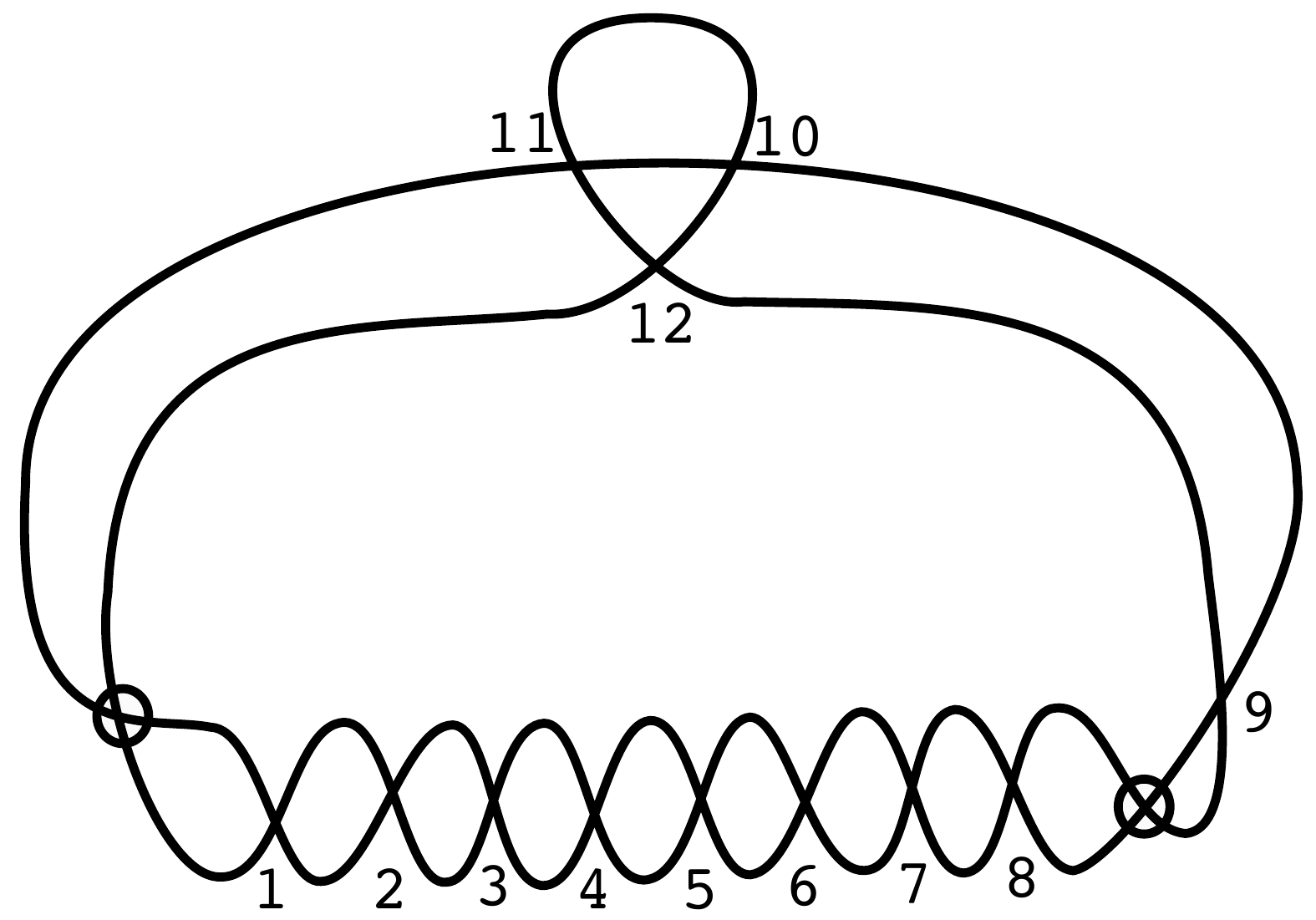}}\hspace{.5in}
\subfloat[][]{
\begin{tabular}{|c|c|}
\hline
$\tr{Y} $& 10\\ \hline
$\cl{Y} $& 8\\ \hline
$\vir{Y}$ & 5\\ \hline
$\kn{Y} $& 4\\ \hline
\end{tabular}}
\hspace{.5in}
\subfloat[][]{\includegraphics[scale=.17]{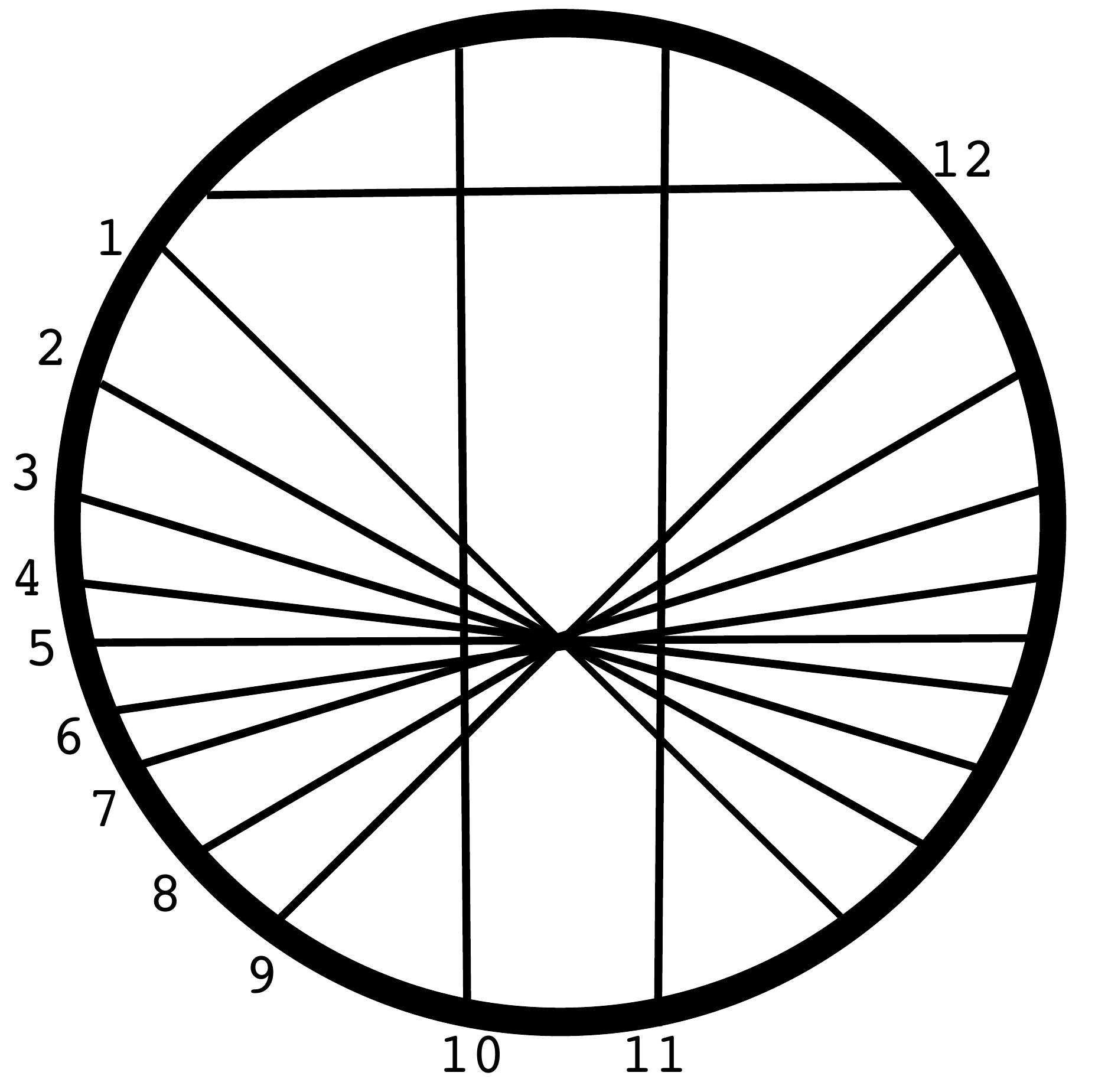}}\hspace{.5in}
\caption{\small{Subfigure (a) shows the virtual shadow $Y$ with classical precrossings labelled. The chord diagram of $Y$ is illustrated in (c), and the the trivializing, classicalizing, virtualizing, and knotting numbers of $Y$ are listed in (b).}}
\label{yay}\label{yaycd}
\end{figure}

We will show that the inequalities in Fact \ref{ingies} are, in fact, sharp. We will also illustrate that classicalizing, trivializing, knotting and virtualizing numbers can be distinct, even for a prime knot's shadow. An example of such a prime shadow is found in Figure \ref{yay}.

\begin{definition}
The \emph{virtual shadow linking number} of two oriented components $L_1$ and $L_2$ in a shadow of a virtual link is the sum of the signs of precrossings between $L_1$ and $L_2$, where the sign of such a precrossing is defined as in Figure \ref{shadowlink}.

\begin{figure}[h]
\centering
\includegraphics[trim=0 500 0 20, width=45mm]{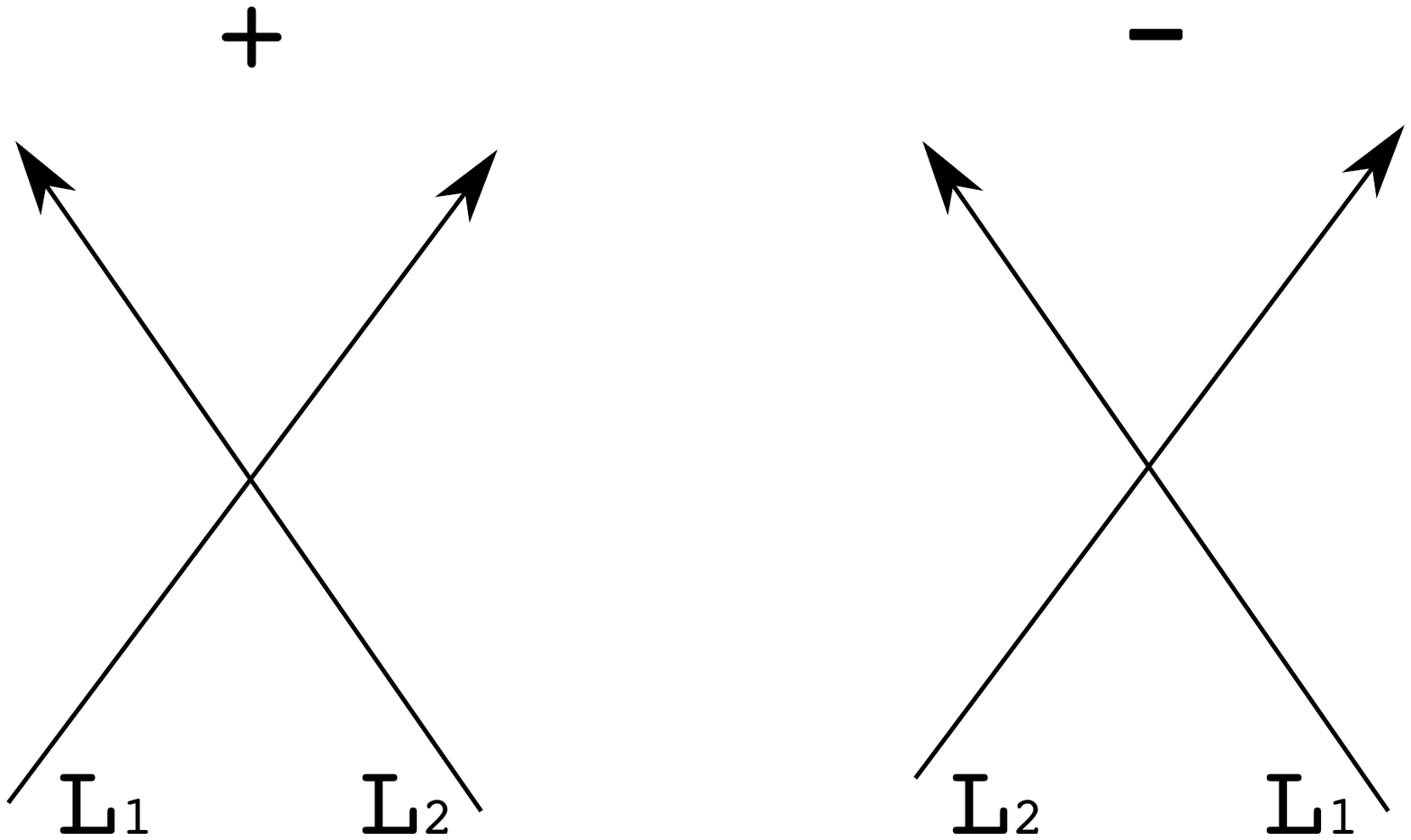}
\caption{Definition of crossing signs for the virtual shadow linking number.}
\label{shadowlink}
\end{figure}
\end{definition}
\begin{definition}
The \emph{intersection index} $ind(c)$ of a crossing or precrossing $c$ in an oriented virtual pseudodiagram is the virtual shadow linking number of the components $L_1$ and $L_2$ created by smoothing at $c$ with the orientation.
\end{definition}

We denote as $C_m$ the set of precrossings and crossings $c$ with intersection index $m$:  $$C_m=\{c:ind(c)=m\}.$$
\begin{definition}[Henrich, \cite{allison}]
The \emph{intersection index polynomial} of a virtual knot $K$ is $$p_t(K) = \sum_{c \in C} w_c \left(t^{|ind(c)|}-1\right),$$ where $C$ is the set of all crossings in $K$ and $w_c$ is the local writhe of crossing $c$.
\end{definition}

As described in \cite{allison}, the intersection index polynomial is an invariant of virtual knots with the property that $p_t(K)=0$ for every classical knot $K$. We now use this polynomial to bound classicalizing and virtualizing numbers.

\begin{proposition}\label{polycling}
For any virtual shadow $S$, $$\cl{S}\ge\sum_{m\ne 0} |C_m|$$ 
\end{proposition}
\begin{proof}
Suppose we resolve fewer than $\sum_{m\ne 0} |C_m|$ precrossings. Then, there exists some precrossing $\hat{c}$ of $S$ such that $\hat{c}\in C_m$ for some $m\ne 0$.  Resolve all other precrossings arbitrarily as classical crossings to obtain $S'$.  The precrossing $\hat{c}$ can be resolved as a positive or negative classical crossing to obtain the knots $S'_+$ and $S'_-$ respectively.  We find $$p_t(S'_+)=t^{|ind(\hat{c})|}-1+\sum_{c \in C\backslash\{ \hat{c}\}} w_c\left(t^{|ind(c)|}-1\right)$$ and
$$p_t(S'_-)=1-t^{|ind(\hat{c})|}+\sum_{c \in C\backslash\{ \hat{c}\}} w_c\left(t^{|ind(c)|}-1\right).$$ 
Since $p_t(S'_+) - p_t(S'_-) = 2(t^{|ind(\hat{c})|}-1)$, and $t^{|ind(\hat{c})|}-1\ne 0$, it follows that at least one of $p_t(S'_+), p_t(S'_-)$ is nonzero.  Thus, our initial resolution of fewer than $\sum_{m\ne 0} |C_m|$ precrossings did not guarantee that the resulting diagram was classical. 
\end{proof}

\begin{lemma}\label{laundry}
For all virtual shadows that have a resolution as a classical knot and for all $m \neq 0$, $|C_m|$ is even.
\end{lemma}
\begin{proof}
Let $S$ be a virtual shadow with the classical resolution $K$. Then, $$p_t(K)=\sum_{c\in C} w_c \left(t^{|ind(c)|}-1\right)=0.$$  In particular, the coefficient of $t^m$ is 0 for all $m\ne 0$. As the coefficient of $t^m$ where $m\ne 0$ is $\sum_{c\in C_m} w_c,$ this implies that $|C_m|$ is even for $m\ne 0$.

Now observe that $ind(c)$ is independent of how $S$ is resolved (when precrossings of $S$ can only be resolved classically).  So for every precrossing $\hat{c}$, $|ind(\hat{c})|$ is constant over every resolution of $S$.  Hence, if $S$ has a classical resolution, then every nonzero $|ind(c)|$ appears an even number of times.
\end{proof}

\begin{proposition}\label{polyvir}
For any virtual shadow $S$ where some precrossing has nonzero intersection index,  
$$\vir{S}\le\min_m \left\{\frac{|C_m|}{2}+1\right\},$$ where $m$ ranges over all positive integers such that there exists a precrossing $\hat{c}$ with $|ind(\hat{c})|=m$.
\end{proposition}
\begin{proof}
Take any positive integer $m$ such that there exists a precrossing $\hat{c}\in C_m$. 
It suffices to show that $\vir{S}\le \frac{|C_m|}{2}+1$ for this particular $m$. By Lemma \ref{laundry}, $\frac{|C_m|}{2} \in \mathbb{Z}$. Designate $\frac{|C_m|}{2}+1$ of the precrossings $c\in C_m$ so that they all have local writhe $+1$. Then, however the other precrossings are resolved, we certainly have that $p_t \neq 0$. Therefore, resolving these $\frac{|C_m|}{2}+1$ precrossings guarantees that any further resolution is non-classical.
\end{proof}

We can now show that the virtual shadow $Y$ shown in Figure \ref{yay} has distinct values for trivializing, classicalizing, virtualizing, and knotting numbers. Note that $Y$ is a prime shadow.

\begin{itemize}

\item $\cl{Y}=8$.

There are eight precrossings, 1--8, with intersection index 2.  All other precrossings have intersection index 0. By Proposition \ref{polycling}, $\cl{Y}\ge 8$.

It is clear that if we resolve the precrossings 1, 3, 5, 7 as positive classical crossings and the precrossings 2, 4, 6, 8 as negative classical crossings, then the resulting pseudodiagram is necessarily classical. Thus, $\cl{Y}=8$.

\item $\tr{Y}=10$.

We see that crossings 1--8 must be resolved to guarantee that $Y$ is classical. Removing these resolved crossings (as specified above) as well as the two flanking virtual crossings using repeated type 2 moves results in a classical knot shadow that has trivializing number 2. Thus, $\tr{Y}=10$.

\item $\vir{Y}=5$.

By Proposition \ref{polyvir}, $\vir{Y} \leq \frac{8}{2}+1=5.$ Suppose only four precrossings are resolved.  Then we can resolve the remaining precrossings from $\{1,2,3,4,5,6,7,8\}$ in a way that allows us to eliminate the virtual crossings through a sequence of Reidemeister II moves. Thus $\vir{Y}=5$.

\item $\kn{Y}=4$.

If we resolve precrossings 9, 10, 11, 12 as positive classical crossings, then the resulting pseudodiagram is necessarily knotted.

It is clear that resolving any 3 precrossings does not guarantee a knotted resolution.  Thus $\kn{Y}=4$.
\end{itemize}

\begin{theorem}
The inequalities in Fact \ref{ingies} are sharp.
\end{theorem}  
\begin{proof}
Consider the virtual shadow of the virtual trefoil shown in Figure \ref{arbitrarytrefcomp}.  It is clear by inspection that $\kn{P}=\vir{P}=\cl{P}=\tr{P}=2.$
\end{proof}

\begin{theorem}
The differences between $\cl{P}$ and $\tr{P}$, and between $\kn{P}$ and $\vir{P}$ can be arbitrarily large.
\end{theorem}
\begin{proof}
Consider a composition of a shadow of a classical trefoil and a 2-braid with one virtual crossing and $2n\ge 4$ classical precrossings, as illustrated in Figure \ref{braidtref}.  
\begin{figure}[h]
\includegraphics[trim=0 450 0 0, width=50mm]{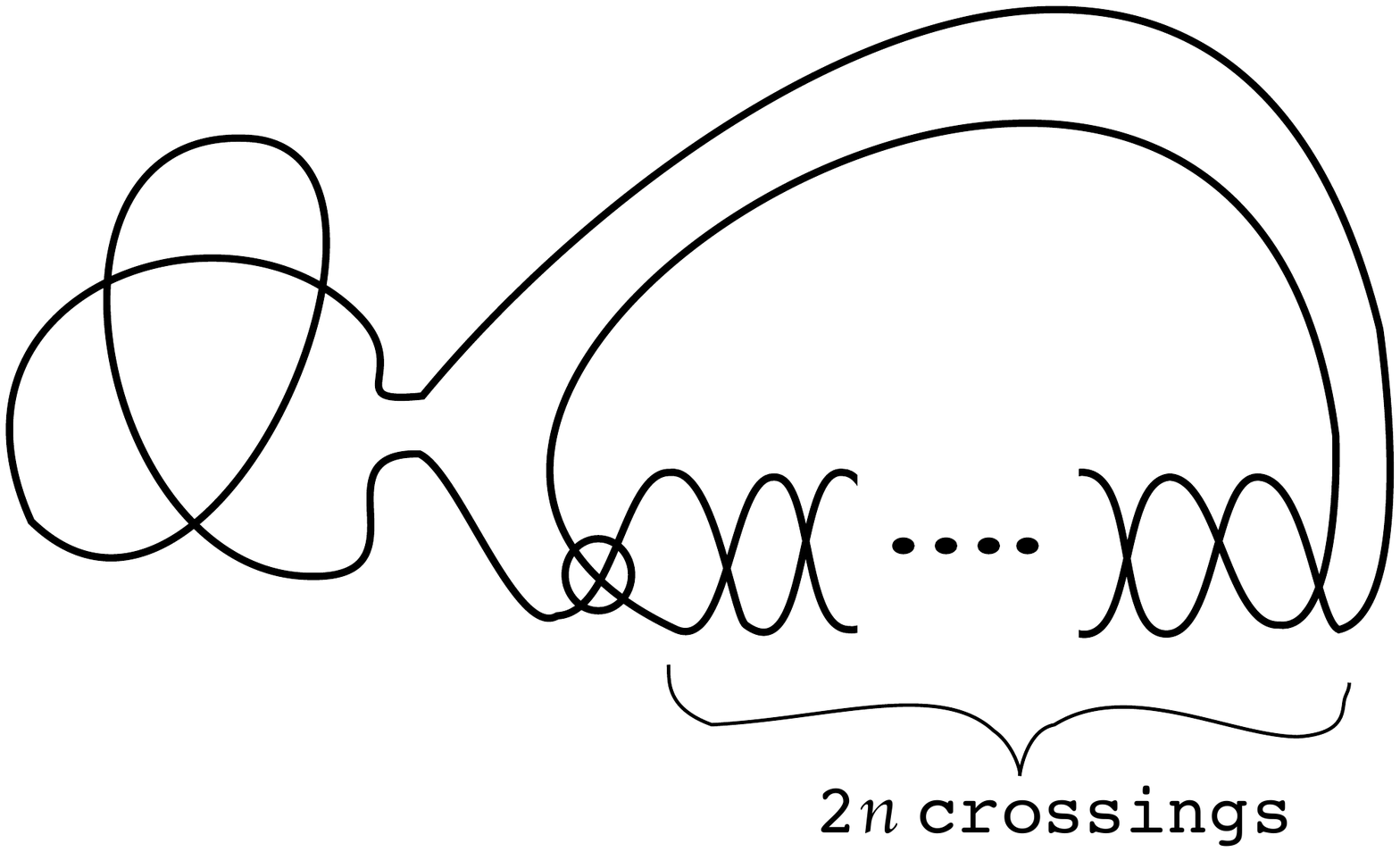}
\caption{\small{The shadow of a 2-braid with one virtual crossing composed with a classical trefoil. This shadow has $\kn{S}=3, \vir{S}=n+1$. }}
\label{braidtref}
\end{figure} 

In~\cite{km}, it is shown that the composition of a trefoil with any virtual knot is non-trivial. Thus, it is clear that $\kn{S}\le 3$ by resolving the crossings in the trefoil shadow.  By Theorem \ref{knottingnotsmall}, $\kn{S} \geq 3$. Hence, $\kn{S} = 3$. 

By Proposition \ref{polyvir}, $\vir{S}\le n+1$, as all of the 2-braid precrossings have intersection index 1, while the trefoil precrossings have intersection index 0.

Suppose that we resolve fewer than $n+1$ precrossings in the knot.  Then in particular, at most $n$ of the $2n$ 2-braid crossings are resolved, and we may resolve the remaining 2-braid precrossings so that the braid unwinds through a sequence of Reidemeister II moves. This yields a classical knot.

Thus, the composition has $\kn{S}=3$, $\vir{S}=n+1$, and $\lim_{n\to\infty}\vir{P}-\kn{P}=\infty$.

The virtual shadow in Figure \ref{arbitrarytrefcomp} shows arbitrarily large differences between $\cl{S}$ and $\tr{S}$.

\begin{figure}[h]
\subfloat[1][A virtual shadow of a virtual trefoil.]{\includegraphics[trim= 0 600 250 0, width=40mm]{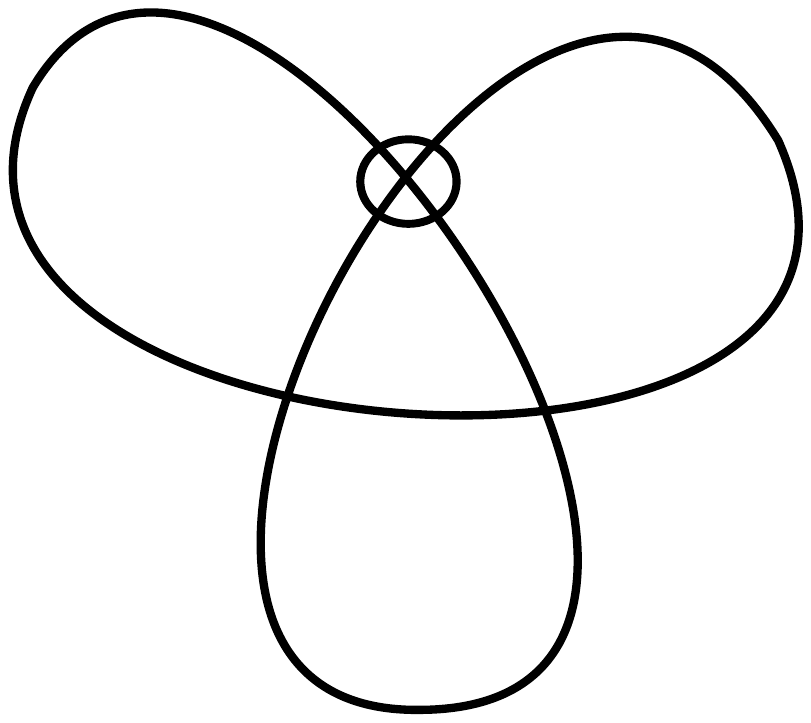}}\hspace{.2in}
\subfloat[2][A shadow composition of an arbitrary number of trefoils with a virtual trefoil shadow.]{\includegraphics[trim= 0 600 0 0, width = 60mm]{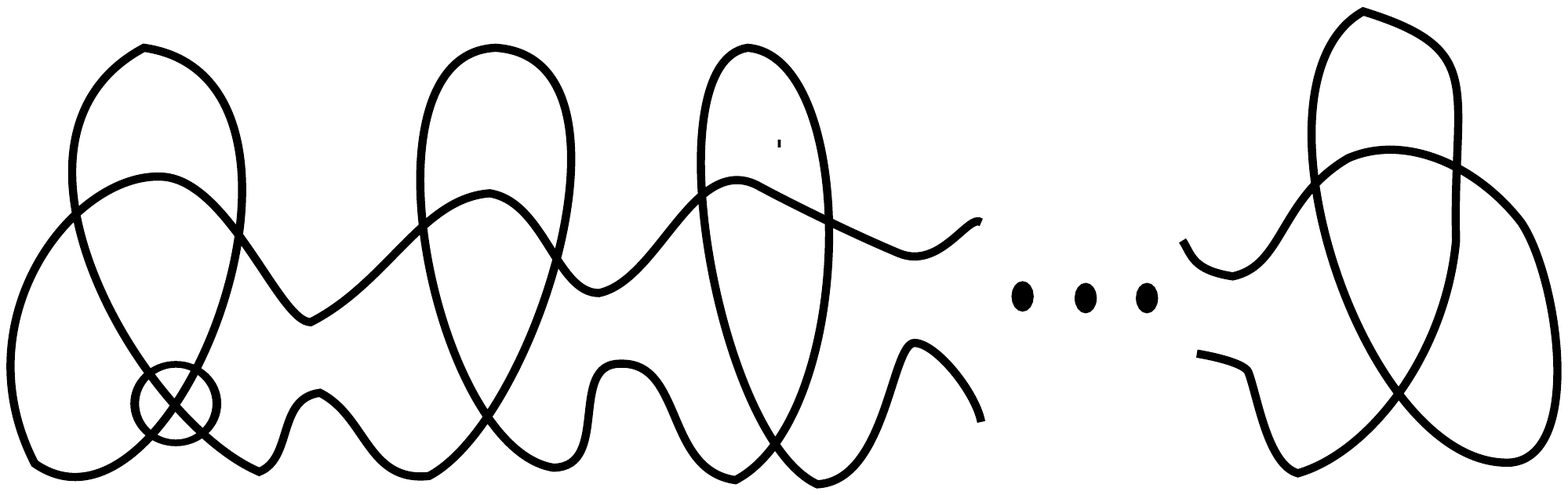}}
\caption{\small{The shadow in (A) has $\kn{P}=\vir{P}=\cl{P}=\tr{P}=2$, while the composition shown in (B) shows the distance between $\cl{P}$ and $\tr{P}$ may be arbitrarily large.}}
\label{arbitrarytrefcomp}
\end{figure}
 The shadow is resolved nontrivially if the virtual trefoil shadow or any one of the $n$ classical trefoil shadows are resolved nontrivially.  Therefore $\tr{S}=2n+2$.

In order to guarantee that $S$ resolves as a classical knot, it suffices to trivialize the virtual trefoil shadow, so $\cl{S}=2$. Hence, $\lim_{n\to \infty}\tr{S}-\cl{S}=\infty$.
\end{proof}

\begin{figure}[h]
\begin{center}
\includegraphics[scale=.5]{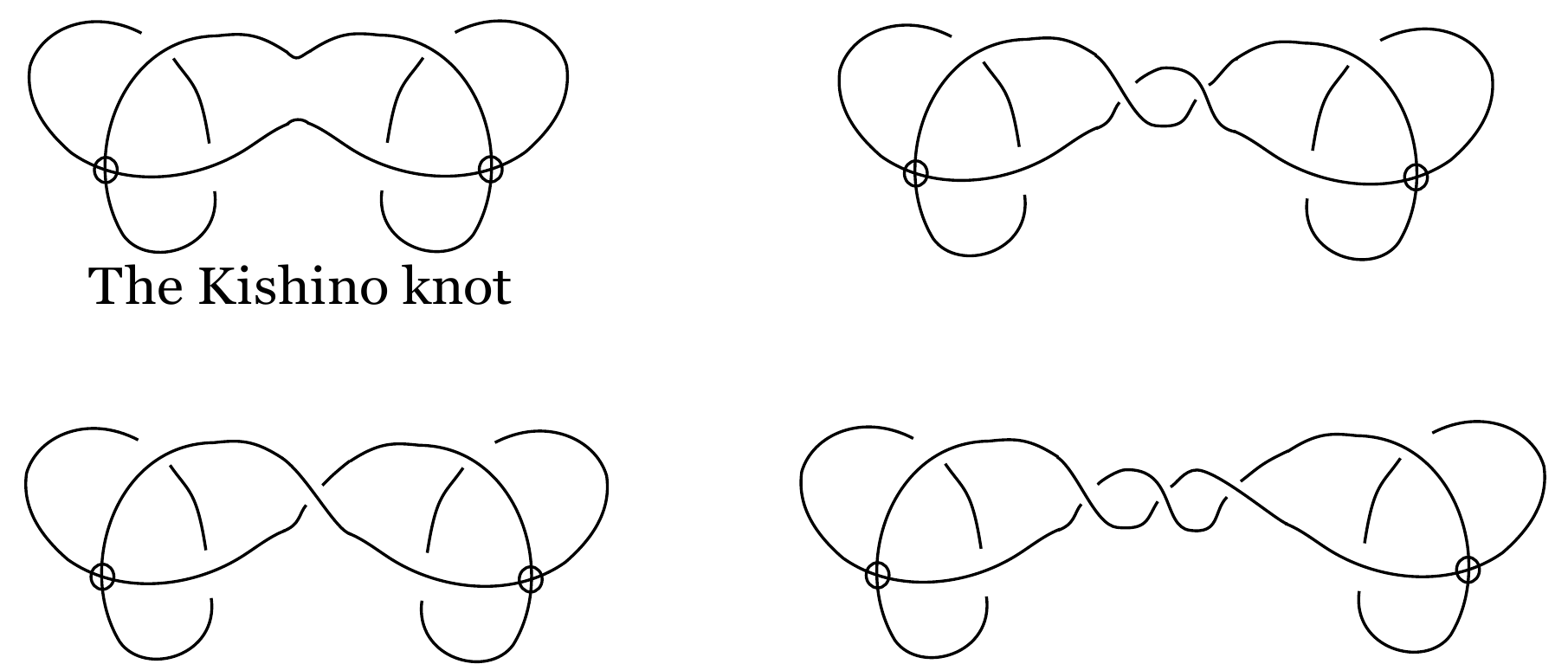}
\end{center}
\caption{\small{The first four members of an infinite family of distinct compositions of two unknots.}}
\label{kishifam}
\end{figure}

Composition is not well-defined for virtual knots. Indeed, there are an infinite number of non-equivalent ways to compose any two knots.  In the virtual realm, composing two unknots may yield a non-trivial and indeed non-classical result, such as the Kishino knot shown in Figure \ref{kishifam}. Here we provide an infinite family of distinct compositions of two unknots. Distinctness can be shown using Kauffman and Dye's arrow polynomial, as defined in~\cite{Arrow}.

\begin{proposition}\label{coolcon}
For any pair of classical knots $K_1$ and $K_2$, there is a shadow $S$ of a composition $K_1\#K_2$ of two virtual diagrams of  $K_1$ and $K_2$ such that $\cl{S}=\infty$.
\end{proposition}
\begin{proof} Consider classical diagrams $D_1$ and $D_2$ of $K_1$ and $K_2$ respectively.  Perform Reidemeister moves as shown in Figure \ref{diagramcomp} (A) to obtain diagrams $D'_1$ and $D'_2$.  The composition of $D'_1$ and $D'_2$ as illustrated in Figure \ref{diagramcomp} (B) yields a diagram that is also a composition of the Kishino knot with the diagrams $D_1$ and $D_2$.
\begin{figure}[h]
\subfloat[][The diagram $D'_1$ of $K_1$]{\includegraphics[trim= 0 50 80 0, width=.3\textwidth]{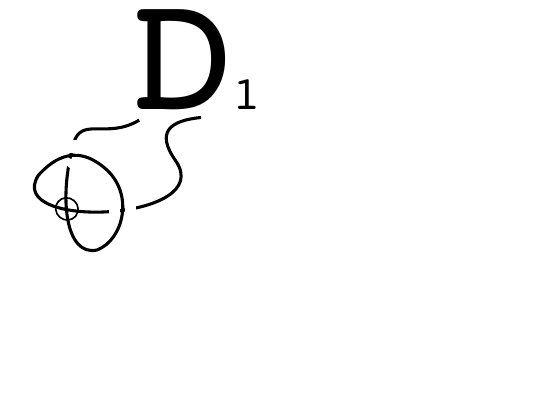}}\hspace{.5in}
\subfloat[][$D'_1\#D'_2$]{\includegraphics[trim= 0 300 0 0, width=.3\textwidth]{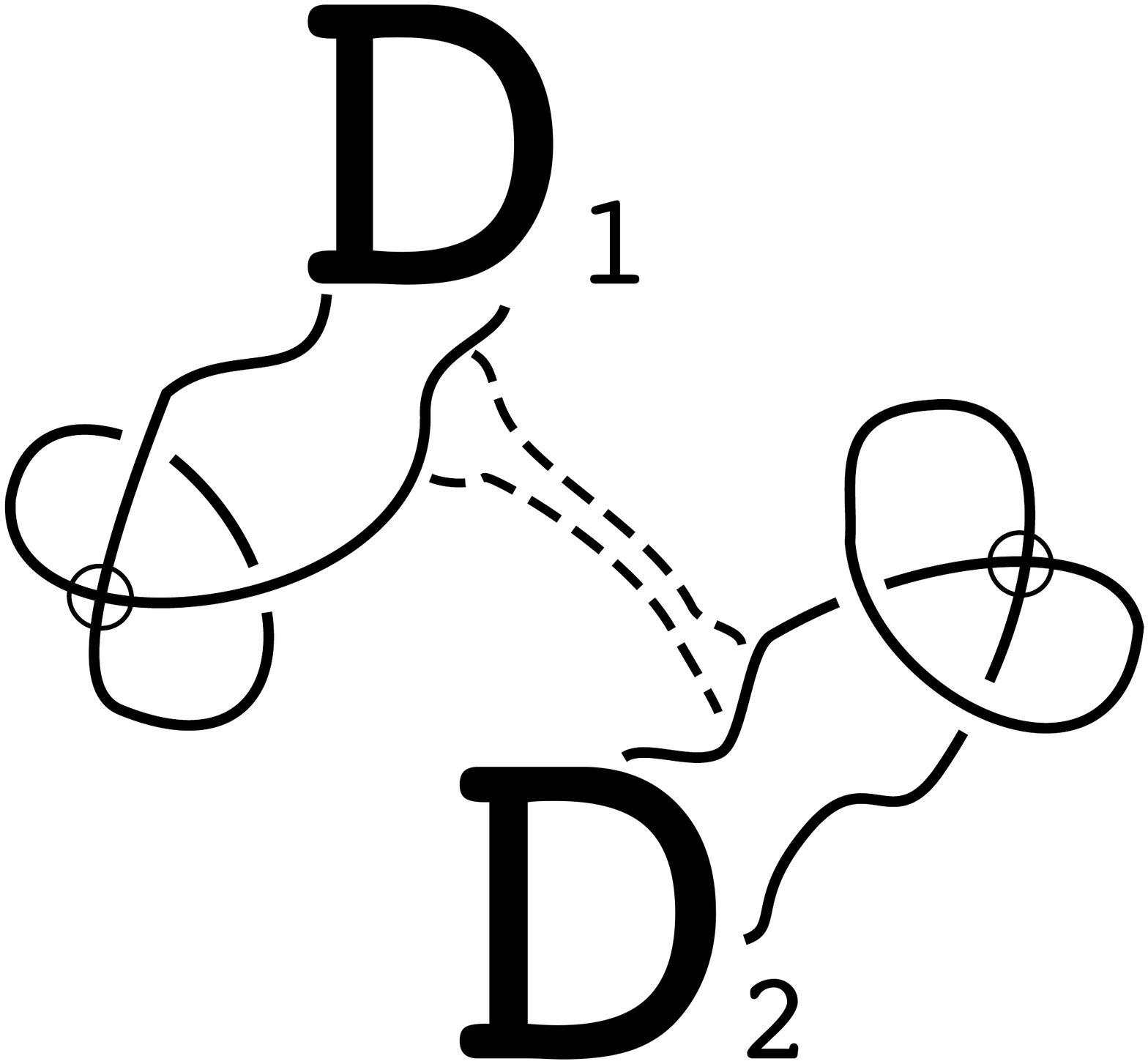} }
\caption{Two stages in the construction of Proposition \ref{coolcon}.}
\label{diagramcomp}
\end{figure}
It is clear from the chord diagrams that, for any two virtual pseudodiagrams $P_1,$ $P_2$, $$\vir{P_1\#P_2}=\min\{\vir{P_1},\vir{P_2}\}.$$
Consider the shadow $S$ of this diagram.  As $S$ is the shadow of the composition of the Kishino knot with two other knots, $\vir{S}\le \vir{\text{Kishino}}=0$.   Thus $\cl{S}=\infty$. 
\end{proof}

We now consider relations between the various trivializing numbers defined in Definitions \ref{tr}, \ref{virtr}, and \ref{ubertr}.
\begin{lemma}\label{tr=vtr}
For any classical shadow $S$, $\virtr{S}\leq\tr{S}$.
\end{lemma}
\begin{proof}
We know that $\tr{S}$ is the least number of chords that must be removed from the chord diagram of $S$ so that the resulting diagram is parallel. Turning the precrossings associated to this deleted set of intersecting chords into virtual crossings will result in a parallel chord diagram. As in the proof of Lemma~\ref{triv=crosses}, we may unknot any such pseudodiagram. In this case, however, we use a sequence of virtual Reidemeister moves to virtually untangle loops before removing each precrossing, in turn, with a classical type I move.



\end{proof}

\begin{lemma}
For any \"{u}ber-virtual shadow $S$, $$\ubtr{S}\leq\virtr{S}\text{  and  }\ubtr{S}\leq\tr{S}.$$
\end{lemma}
\begin{proof}
If we can trivialize $S$ by designating the crossings in a crossing set $\mathscr{T}$ as virtual or by designating the crossings in $\mathscr{T}$ as classical, then we can certainly trivialize $S$ by setting the crossings of $\mathscr{T}$ as virtual \emph{or} classical.  

\end{proof}

\begin{theorem} For any \"{u}ber-virtual shadow $S$, $\ubtr{S}=\tr{S}$. Hence, $$\ubtr{S}=\virtr{S}=\tr{S}.$$
\end{theorem}
\begin{proof} In Section~\ref{Kauffman's $J$-invariant}, we will see that any pseudodiagram with a chord diagram that has intersecting precrossings can be resolved non-classically. Hence, $\ubtr{S}=\tr{S}$. The second equality is immediate from the previous lemmas.
\end{proof}

\begin{proposition}
For any virtual shadow $S$, $\tr{S} \neq 1$.
\end{proposition} 
\begin{proof}
Let $S$ be a virtual shadow with $\tr{S} \leq 1$.  Then, there is a trivial virtual pseudodiagram $T$ obtained from $S$ by a classical resolution of some precrossing $p$.  Without loss of generality, assume $p$ was resolved as a positive classical crossing. 

Now, let $T'$ be the mirror image of $T$, where we assume the mirror image of a precrossing is still a precrossing, and the mirror image of a virtual crossing is a virtual crossing. Since $T$ is trivial, we have by symmetry that $T'$ is also trivial. However, $T'$ is precisely $S$ with $p$ resolved as a negative classical crossing. Therefore, $S$ is trivial regardless of how we resolve $p$ classically, and so $\tr{S} = 0$.
\end{proof}

\section{Kauffman's $J$-invariant}
\label{Kauffman's $J$-invariant}
In this section we will make use of Kauffman's $J$-invariant for virtual knots, defined in \cite{Kauffman}.  Kauffman's $J$-invariant for a virtual knot $K$ is defined as $$J(K)=\sum_{c\in Odd(K)}w_c,$$ where $w_c$ is the local writhe of crossing $c$, which corresponds to a chord in $Odd(K)$. 

We say that $c\in Odd(K)$ if an arc of $K$ with both endpoints at $c$ passes through an odd number of classical crossings.  Equivalently, the chord corresponding to $c$ in the chord diagram of $K$ intersects an odd number of chords.

For all classical diagrams $D$, $Odd(D)=\emptyset$, and therefore $J(D)=0$.  Hence, if $J(D)\ne 0$, then $D$ is necessarily non-classical.

\begin{proposition}\label{odd}
For any virtual pseudodiagram $P$ containing a precrossing in $Odd(P)$, $P$ can be resolved non-classically (and therefore, non-trivially).
Moreover, $\cl{P}\ge |Odd(P)\cap C_p|$, where $C_p$ is the set of precrossings in $P$.
\end{proposition} 
\begin{proof}
Let $$W=\sum_{c\in Odd(P)\cap \overline{C}_p}w_c.$$  We can resolve all precrossings in $Odd(P)$ with local writhe $+1$ such that $J(P)=W+|Odd(P)\cap C_p|$, or with local writhe $-1$ such that $J(P)=W-|Odd(P)\cap C_p|.$  If there exists a precrossing in $Odd(P)$, then at least one of these resolutions has nontrivial $J$, and therefore is non-classical.
\end{proof}

\begin{lemma}
For any virtual pseudodiagram $P$, $|Odd(P)|$ is even.
\end{lemma}

\begin{proof}
This result is a straightforward consequence of graph theory.
\end{proof}






\begin{theorem}
For any virtual pseudodiagram $P$ with $|Odd(P)| \neq 0$, $$\vir{P} \le \frac{|Odd(P)|}{2} + 1.$$
\end{theorem}
\begin{proof} 
In order to make sure that any further resolution of $P$ is necessarily non-classical, it suffices to resolve some precrossings of $P$ such that any further resolution to a diagram $D$ has $J(D) \neq 0$.  Thus, it suffices to resolve $\frac{|Odd(P)|}{2} + 1$ chords from $Odd(P)$ to have local writhe $+1$ to obtain $P'$. Regardless of how the other chords are resolved, $\sum_{c \in Odd(P')} w_c > 0$, and $P'$ is necessarily non-classical.
\end{proof}

To more generally analyze the trivializing number of a virtual pseudodiagram, we might begin by asking whether the proof of Lemma~\ref{+nontriv} holds for virtual pseudodiagrams. In general, the answer is negative for \emph{compact} virtual pseudodiagrams, that is, for the sorts of virtual pseudodiagrams we've been discussing thus far. This is because the Vassiliev invariant used in the proof is no longer base-point independent, and hence is not an invariant for virtual knots. If we were to consider \emph{long} virtual psudodiagrams, or based virtual pseudodiagrams, the result would extend readily. Fortunately, we may repair the damage by using Kauffman's $J$-invariant and several recent results due to Manturov.

\begin{lemma}\label{ManturovProof} Any virtual pseudodiagram $P$ with a chord diagram containing intersecting prechords can be resolved nontrivially.\end{lemma}

\begin{proof}If every crossing in a virtual pseudodiagram is even, then the Vassiliev invariant used in Lemma~\ref{+nontriv} is invariant mod 2. This relies on the following result from~\cite{freeknot}. If two equivalent virtual knot diagrams have even crossings, then the virtual knot diagrams are related by a sequence of Reidemeister--type moves involving only even crossings. Thus, Lemma~\ref{+nontriv} holds for virtual pseudodiagrams containing only even crossings.

Now, let us first resolve all precrossings in $P$ so that exactly two precrossings corresponding to intersecting prechords in the chord diagram remain. Assume that at least one of these precrossings is odd. Then, by Proposition~\ref{odd}, we can choose crossing information for the odd precrossing(s) so that the resulting virtual knot is non-classical. So let us turn to the remaining case where our virtual pseudodiagram contains odd crossings, but the two precrossings are even.

In~\cite{freeknot}, Manturov introduces a functorial mapping, $f$, on the set of virtual knots that sends a virtual knot diagram $K$ to the virtual knot diagram $f(K)$ where all odd classical crossings are made virtual. We note that if $K$ is trivial, then so is $f(K)$. So let us examine $f(P)$ where $P$ is a virtual pseudodiagram with two even precrossings that intersect as chords in the chord diagram of $P$. Since our precrossings are both even in $P$, they persist in $f(P)$. They may, however, be odd in $f(P)$. If at least one of the precrossings becomes odd, then we find (as above) that $f(P)$ can be resolved nontrivially, so $P$ can be resolved nontrivially. On the other hand, if both precrossings remain even, we repeat the process of applying the functor $f$ until we have a virtual pseudodiagram with all even crossings and our precrossings persist. In this case, we showed above that the pseudodiagram $f(f(\cdots f(P)\cdots))$ can be resolved nontrivially. Hence, $P$ can be resolved nontrivially. 
\end{proof}

\begin{lemma}\label{+nonclung}
Let $P$ be an \"{u}ber-virtual pseudodiagram with exactly two precrossings $a$ and $b$ such that the chords $a$ and $b$ cross in the chord diagram of $P$.  Then $P$ has a resolution which is non-classical.
\end{lemma}
\begin{proof}
Recall that $J(K)=\sum_{c\in Odd(K)}w_c$ for any knot $K$, and $J=0$ for all classical knots.  

Suppose that either $a$ or $b$ is in $Odd(\alpha)$ for the chord diagram $\alpha$ containing all classical and precrossings of $P$.  Without loss of generality, assume $a\in Odd(P)$. Then for any knot $K$ which is resolved from $P$ with $a$ and $b$ resolved classically, $$J(K)=w_a + \sum_{c\in Odd(\alpha)\backslash\{ a \} } w_c.$$
Therefore, we may resolve $a$ and $b$ as classical crossings with the appropriate local writhes, obtaining a diagram $K$ with $J(K)\ne 0$. Hence, $K$ is non-classical.

Suppose now that both $a$ and $b$ are not in $Odd(\alpha)$.  Resolve $b$ virtually. Then, $a$ is in $Odd(\alpha')$ for the new chord diagram $\alpha'$.  Again, $$J(K)=w_a + \sum_{c\in Odd(\alpha') \backslash \{ a \} }w_c$$ where $K$ is a resolution of $P$ with $b$ resolved as a virtual crossing and $a$ is resolved as a classical crossing. We may resolve $a$ with the appropriate local writhe such that $J(K)\ne 0$.  Hence, $K$ is a non-classical resolution of $P$.

In either case there is a resolution of $P$ which is non-classical. Therefore any \"{u}ber-virtual pseudodiagram with exactly two precrossings that cross in the chord diagram can be resolved non-classically.
\end{proof}
\begin{corollary}
\label{cor1}
If an \"{u}ber-virtual pseudodiagram has precrossings that intersect in the chord diagram, then that \"{u}ber-virtual pseudodiagram can be resolved non-classically.
\end{corollary}

\begin{corollary}\label{vircl=virtr}
For any \"{u}ber-virtual shadow $S$, $\vircl{S}=\virtr{S}$.
\end{corollary}

\section{Unknotting Numbers and Genus}
\label{Unknotting Numbers and Genus}

\begin{definition}
The \emph{unknotting number} $\un{D}$ of a diagram is the minimum number of classical crossings whose local writhes must be switched so that the resulting diagram is a diagram of the unknot.

The unknotting number $\un{K}$ of a knot $K$ is $\min_{D}{\un{D}},$ where $D$ ranges over all diagrams of $K$.  
\end{definition}

More information about the unknotting number can be found in \cite{knotbook}. We define the trivializing number of a knot $K$ to be $\min_{S_K} \tr{S_K}$, where $S_K$ ranges over all shadows of conformations of $K$.  

It is known that the unknotting number cannot always be realized in a minimum crossing projection of a knot. The following is a related open question.

\begin{question}
Is the trivializing number of every knot $K$ realized in a minimum crossing projection of $K$?
\end{question}

Theoretical lower bounds on unknotting number are known. However, upper bounds are generally determined by trial and error. We believe that the following theorem is the first non-trivial upper bound to be developed. 

\begin{theorem}
For any knot $K$, $$\un{K}\le \frac{\tr{K}}{2}.$$
\label{unknottingbound}
\end{theorem}
\begin{proof}
Let $S$ be an oriented shadow of a diagram $D$ of $K$ which realizes the unknotting number.  It suffices to show that $\un{K}\le \frac{\tr{S}}{2}$.  Let $C$ be a list of $\tr{S}$ precrossings $c_1,c_2,\ldots, c_{\tr{S}}$ that can be resolved to trivialize $S$.
We denote the trivializing resolutions of $C$ as $w_1,w_2,\ldots, w_{\tr{S}}$, where $w_i\in \{+1,-1\}$ is the resolved local writhe of $c_i$.  As a knot is nontrivial if and only if its mirror image is nontrivial, $C$ also trivializes $S$ when every $c_i$ is resolved with writhe $-w_i$.

Look at the set $C$ in $D$.  In order to unknot $D$ it suffices to change crossings of $C$ such that 
$c_1,c_2,\ldots, c_{\tr{S}}$ either have local writhes $w_1,w_2,\ldots,w_{\tr{S}}$ or $-w_1,-w_2,\ldots, -w_{\tr{S}}$.  Observe that at most half of the crossings in $C$ must be changed to effect this, as if it requires changing $n$ crossings to obtain local writhes $w_1,w_2,\ldots,w_{\tr{S}}$, then it takes $\tr{S} - n$ changes to    obtain local writhes $-w_1,-w_2,\ldots, -w_{\tr{S}}$. Thus, $\un{K}\le \frac{\tr{S}}{2}$.
\end{proof}

The following invariant related to the unknotting number was introduced by Goussarov, Polyak and Viro in \cite{PolyakViro} and Fleming and Mellor in \cite{introsgt}.

\begin{definition}
The \emph{virtual unknotting number} $\viru{D}$ of a diagram is the minimum number of classical crossings that must be made virtual so that the resulting diagram is a diagram of the unknot.

The virtual unknotting number $\viru{K}$ of a knot $K$ is $$\min_{D}{\viru{D}},$$ where $D$ ranges over all diagrams of $K$.  
\end{definition}

\begin{theorem}
For any knot $K,$ $$\viru{K}\le \tr{K}.$$
\end{theorem}
\begin{proof}
Let $S$ be any shadow of any conformation of $K$.  By Theorem \ref{tr=vtr},  $\tr{S}=\virtr{S}$, so it suffices to show that $\viru{K}\le \virtr{S}$.

The virtual unknotting number of $K$ is the minimum over all diagrams of $K$ of the minimum number of crossings which must be made virtual so that $K$ becomes the unknot.  The virtual trivializing number of $S$ is the minimum number of precrossings in $S$ that must be resolved as virtual to guarantee that any diagram obtained by resolving the remaining precrossings is the unknot.  The result follows.\end{proof}

The similarity of virtual and classical bounds prompts us to consider the relationship between virtual unknotting and unknotting number, a relationship first considered by Fleming and Mellor in \cite{introsgt}.

\begin{theorem}\label{noel}
For any virtual (or classical) knot, $K$, $\viru{K}\leq 2 \un{K}$. 
\end{theorem}
\begin{proof}
Let $D$ be a diagram of $K$ that realizes the classical unknotting number.  At every crossing $c$ whose writhe must be switched to unknot $D$, we first perform a classical Reidemeister II move, then designate two adjacent crossings as virtual to enable a virtual Reidemeister II move, leaving a single crossing with opposite writhe of $c$.  This sequence of moves is illustrated in Figure \ref{switchcross}.

\begin{figure}[h]
\includegraphics[width=0.4\textwidth]{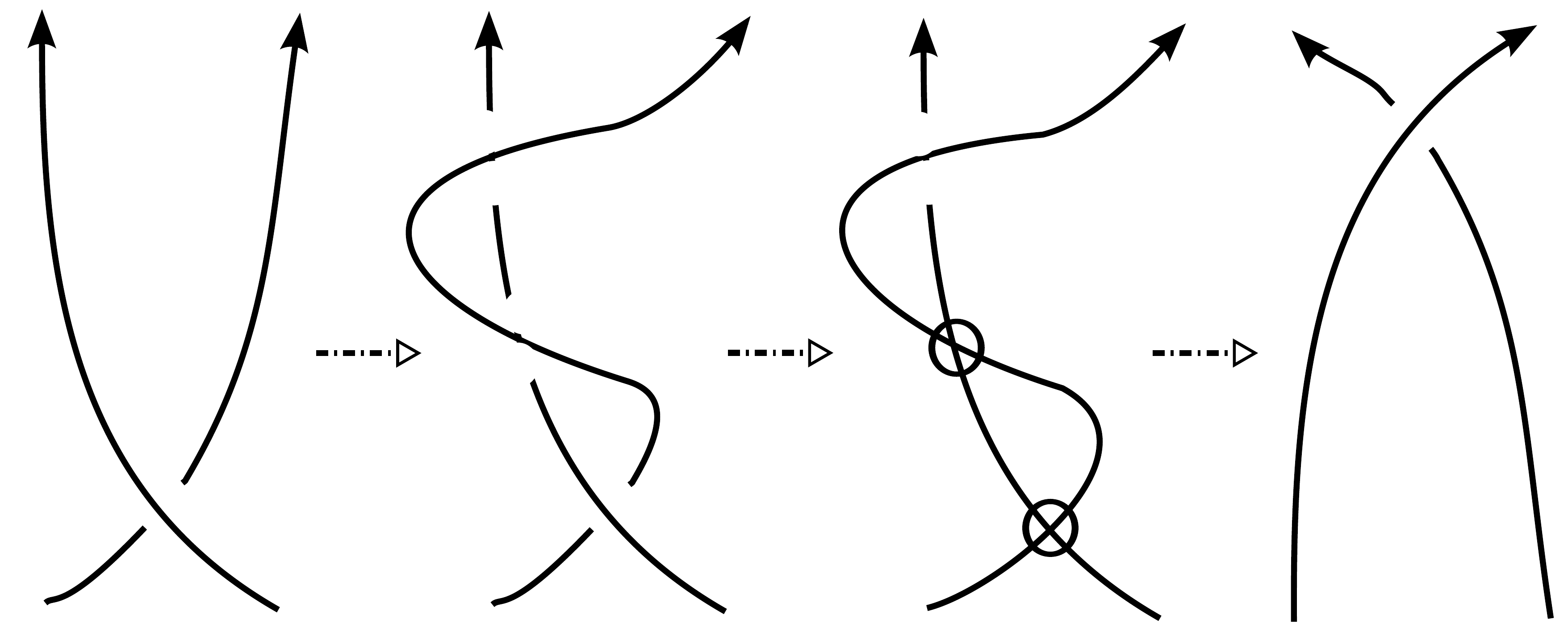}
\caption{\small{This sequence of Reidemeister moves and virtual crossing designations switches the sign of a crossing.}}
\label{switchcross}
\end{figure}

After performing a virtual Reidemeister II move, we obtain the same diagram as would have been obtained by merely switching crossing $c$ in $D$.  Carrying out this process at every crossing that needs to be changed to unknot $D$ results in the same diagram as would be obtained by switching all those crossings.  Thus, $\viru{K} \leq 2 \un{K}$.
\end{proof}

\begin{conjecture}\label{2u=uv}
For all classical knots $K$, $\viru{K} = 2 \un{K}$.
\end{conjecture} 

\begin{figure}[h]
\includegraphics[trim = 0 500 0 0, width=.8\textwidth]{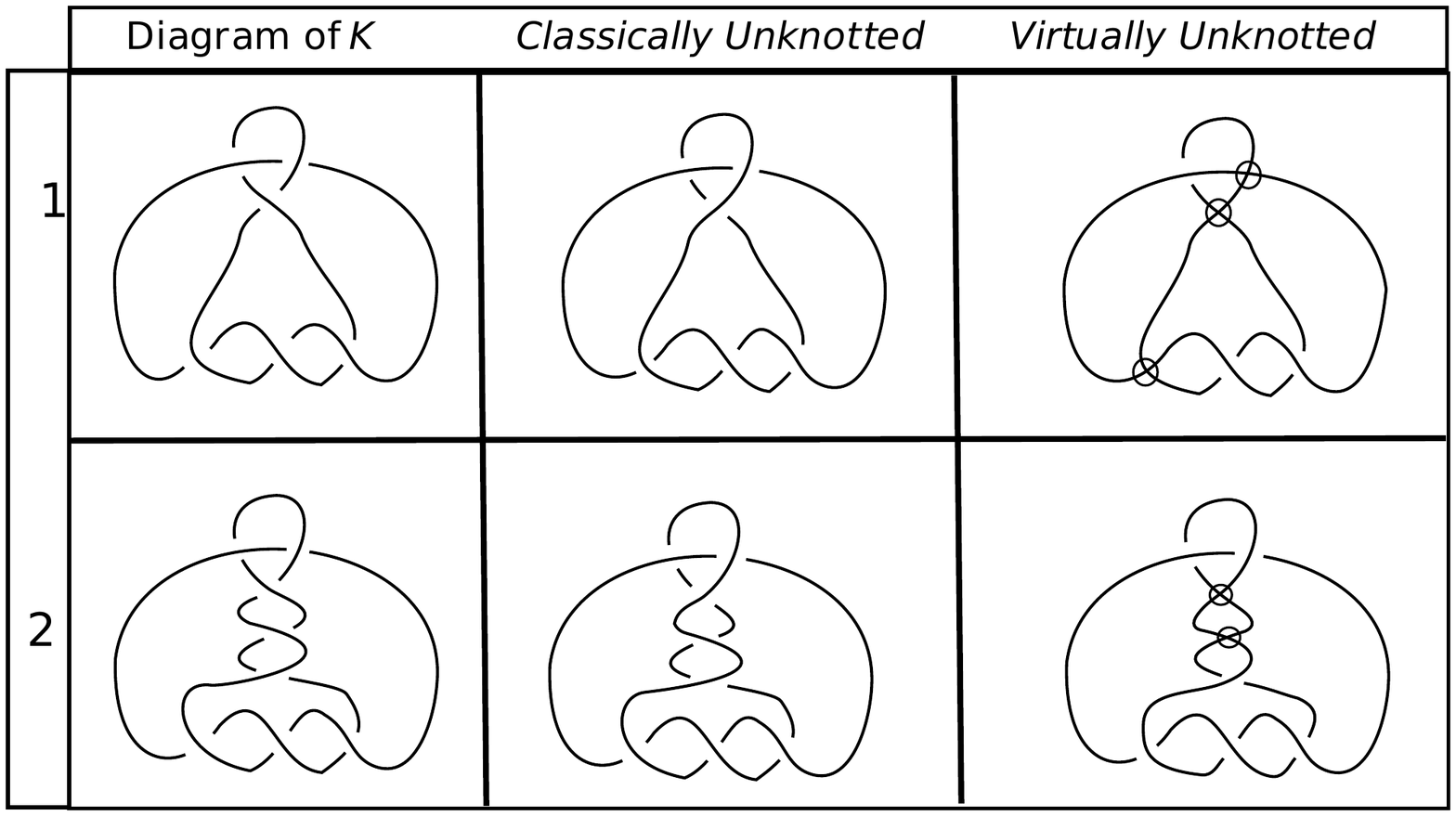}
\caption{\small{Row 1 shows a minimum crossing diagram, with realizations of minimum unknotting and virtual unknotting crossing changes for this diagram.
Row 2 shows a non-minimum crossing projection of the same knot that realizes both the unknotting and virtual unknotting number for the knot.}}
\label{notvunknot}
\end{figure}

There are examples of minimum crossing classical projections that do not realize the unknotting number for a given knot. See for example, \cite{bleiler} and \cite{nakanishi} .  Similarly, there are minimum crossing projections that do not realize the virtual unknotting number.  In Figure \ref{notvunknot}, we exhibit a minimum crossing projection that realizes the unknotting number, but not the virtual unknotting number. Currently there are no examples of diagrams that realize the virtual unknotting number but do not also realize the unknotting number.

Resolving Conjecture \ref{2u=uv} positively would be extremely useful. In particular, the following theorem would not be conditional. 

\begin{theorem}
Assuming either that the unknot is the only classical knot with trivial Jones polynomial, or that no classical knot has virtual unknotting number 1, it follows that the trivializing number for any virtual shadow with precisely one virtual crossing is even.
\end{theorem}
\begin{proof}
Let $S_v$ be a virtual shadow with precisely one virtual crossing $v$, and $S$ be the corresponding shadow with $v$ as a precrossing.
Let $\mathscr{T}=\{c_1,\ldots, c_n\}$ be a minimum trivializing set of precrossings for $S_v$.  We show that $\mathscr{T}$ is a basic trivializing set for $S$, and therefore, by Theorem \ref{H2}, $|T|$ is even.

Let $\mathscr{T}^\star$ denote a resolution of the precrossings of $\mathscr{T}$ such that $S_v$ is trivialized.  Suppose by way of contradiction that resolving $\mathscr{T}$ as $\mathscr{T}^\star$ does not trivialize $S$.  Then there is a knotted resolution $K$ of $S$ where $\mathscr{T}$ is resolved as $\mathscr{T}^\star$. 

If Conjecture \ref{2u=uv} is true, then since $\un{K}\ge 1$, it follows that $\viru{K}\ge 2$. However, $\viru{K}=1$, as changing $v$ to a virtual crossing necessarily produces the unknot.   This would be our desired contradiction.

Alternatively, suppose that $K$ has nontrivial Jones polynomial. Fleming and Mellor proved in \cite{intrinsic} that for diagrams of classical knots with nontrivial Jones polynomial, if one crossing is made virtual, then the Jones polynomial is still nontrivial.  Therefore, for any classical knot with nontrivial Jones polynomial, if one crossing is made virtual, the result is a nontrivial knot.  

Hence, if $v$ is made virtual, $K$ must remain knotted.  However, this contradicts the fact that $S_v$ is trivial when $\mathscr{T}$ is resolved as $\mathscr{T}^\star$. Therefore, assuming that either of the conjectures holds, we conclude that $\mathscr{T}$ trivializes $S$.

\begin{figure}[h]
\centering
\includegraphics[width=.29\textwidth]{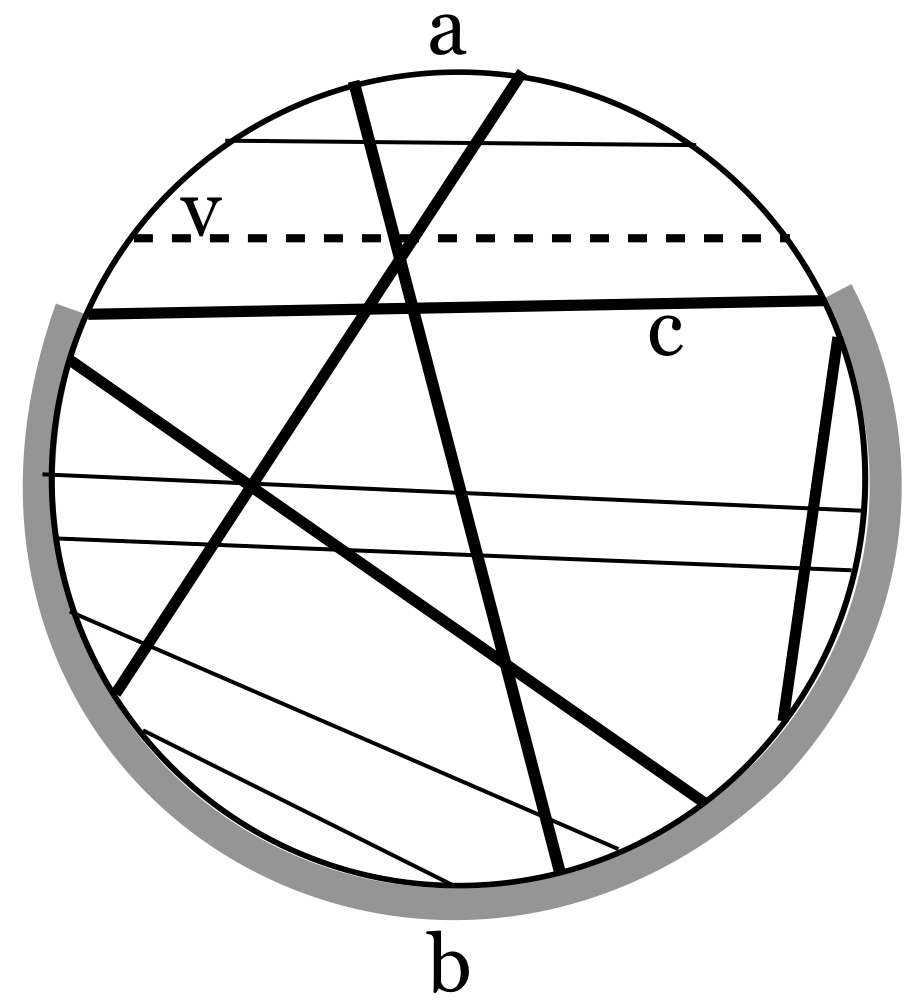}
\caption{\small{The chord diagram of a shadow, with a trivializing set $\mathscr{T}$ in bold, and the crossing $v$ with a dashed line.  By our claim, $\mathscr{T}$ is not basic, so there exists a chord ${\bf c}\in \mathscr{T}$ that intersects only elements of the trivializing set $\mathscr{T}$. Thus $c$ separates the core circle into the black (a) and grey (b) arcs.  Only a crossing in $\mathscr{T}$ may have chord endpoints in both arcs.}}
\label{Trivset}
\end{figure}

Now we show that $\mathscr{T}$ is a \emph{basic} trivializing set of precrossings for $S$.  Suppose by way of contradiction that a proper subset of $\mathscr{T}$ trivializes $S$. (We will show that $\mathscr{T}$ cannot be a basic trivializing set of $S_v$.) By assumption, there is some $c\in \mathscr{T}$ which need not be resolved in $S$ in order to guarantee triviality.  The chord for $c$ in the chord diagram of $S$ does not cross any chords from $S\backslash \mathscr{T}$.  It follows that the chord for $c$ describes two arcs in the core circle where each precrossing chord has both endpoints on a single arc. In particular, as the precrossing $v$ is not an element of $\mathscr{T}$, the endpoints of the chord for $v$ both lie on a single one of these arcs.  We call the arc containing these endpoints $a$, and the other arc $b$.  This configuration is illustrated in Figure \ref{Trivset}.

We claim that there is a resolution of the precrossings $\mathscr{T} \backslash\{c\}$ so that $S_v$ is trivial.  The arcs $a$ and $b$ define two loops $A$ and $B$ based at $c$ that intersect only at crossings in $\mathscr{T}$.  As no endpoint of the chord for $v$ appears along arc $b$, the crossing $v$ does not appear on the corresponding loop $B$.  We say that a point $p$ lies exterior to $B$ if there is a path from $p$ to a point $p'$ arbitrarily far away from the shadow, where this path does not intersect $B$.

Suppose that $v$ lies exterior to $B$.  Then we resolve all crossing between $B$ and $A$ such that $B$ lies completely above $A$.  Because $v$ lies exterior to $B$, it is then possible to perform a sequence of classical Reidemeister moves so that the only crossing between $B$ and $A$ is $c$.  It is then clear that resolving $c$ does not affect the triviality of $S$.

\begin{figure}[h]
\subfloat[1][Concentric Kinks]{\includegraphics[trim=0 250 0 0, width= 30mm]{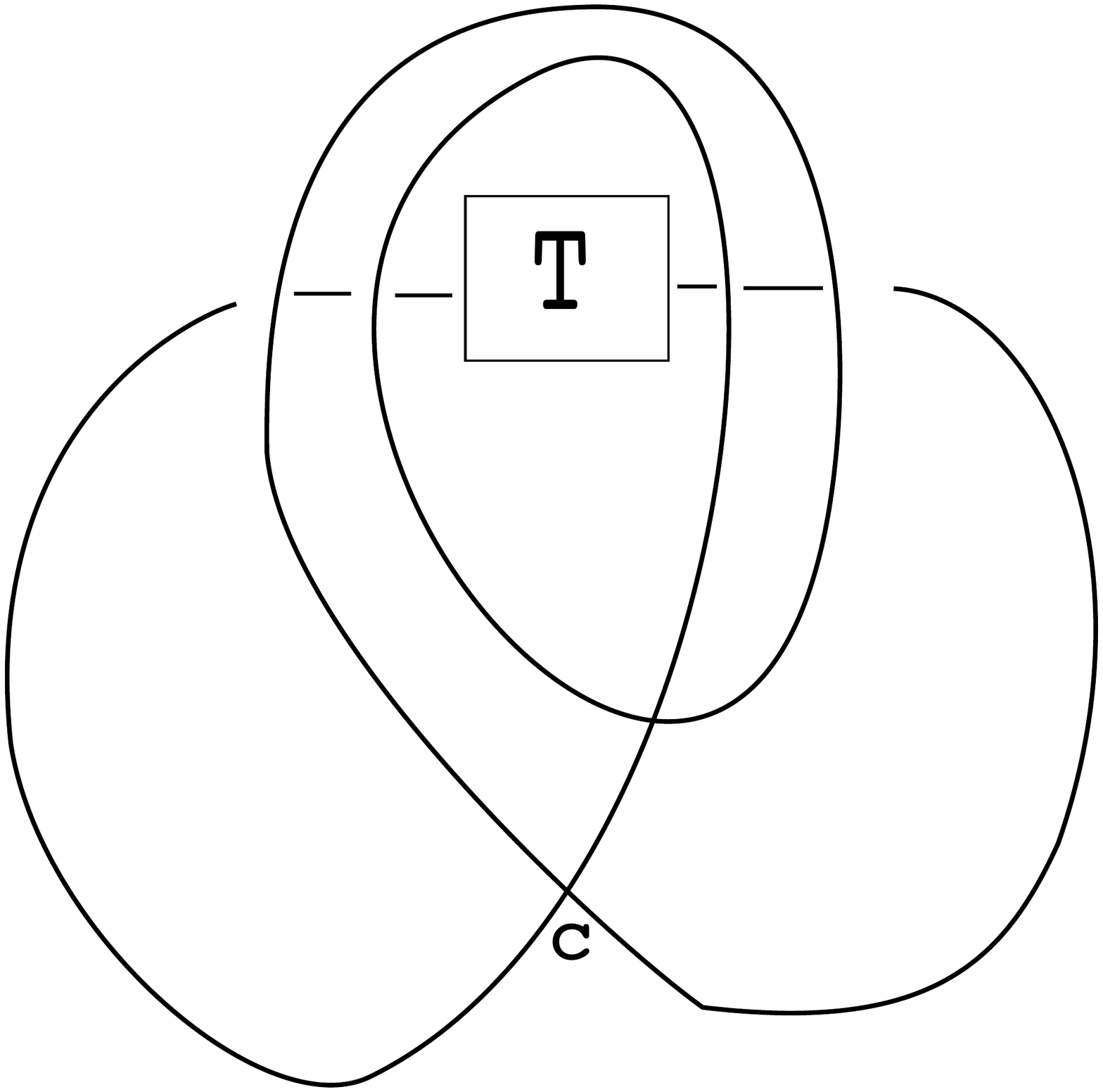}}
\subfloat[2][Unpeeling 1]{\includegraphics[trim=0 250 0 0, width= 30mm]{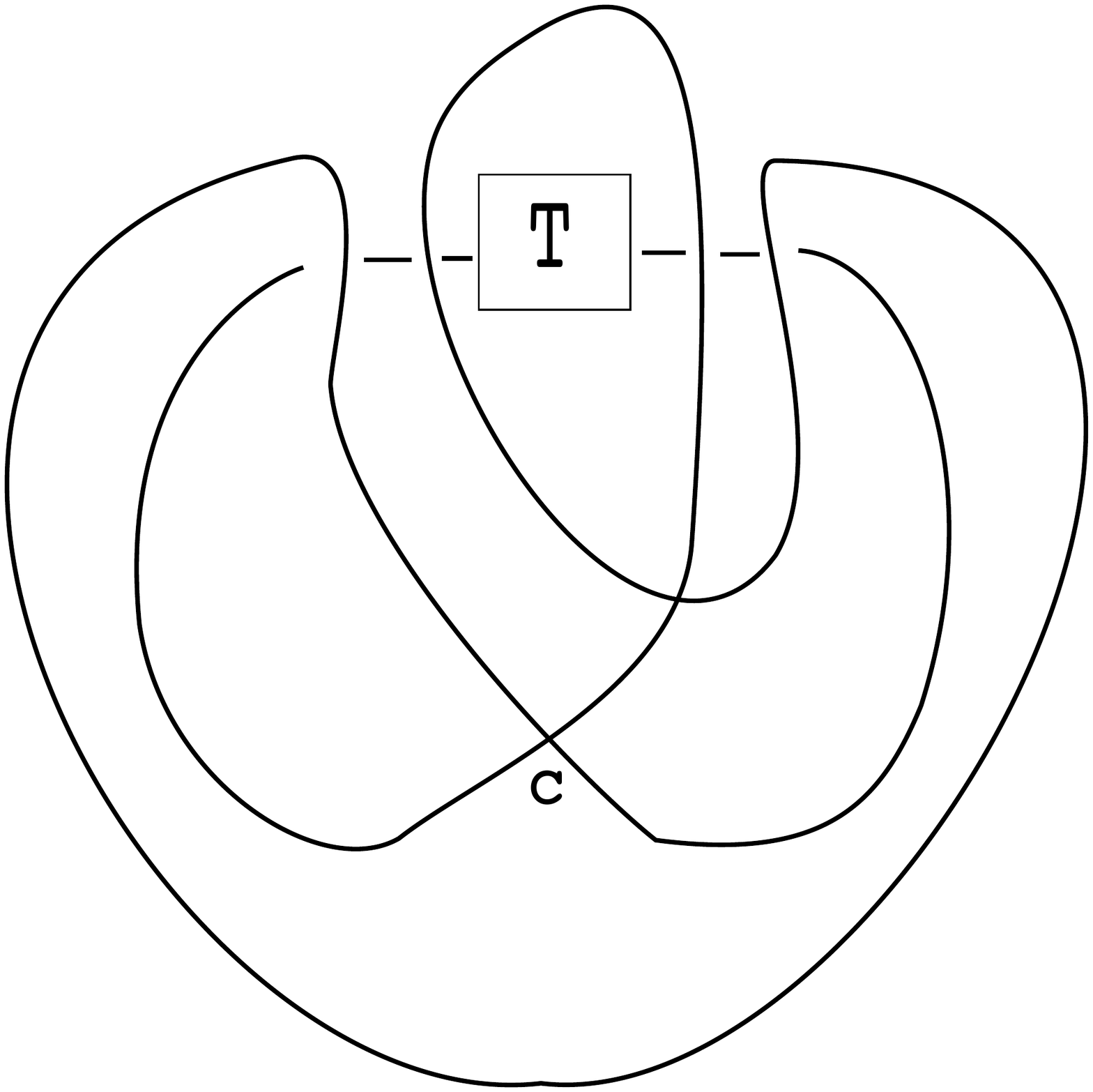}}
\subfloat[3][Unpeeling 2]{\includegraphics[trim=0 250 0 0, width= 30mm]{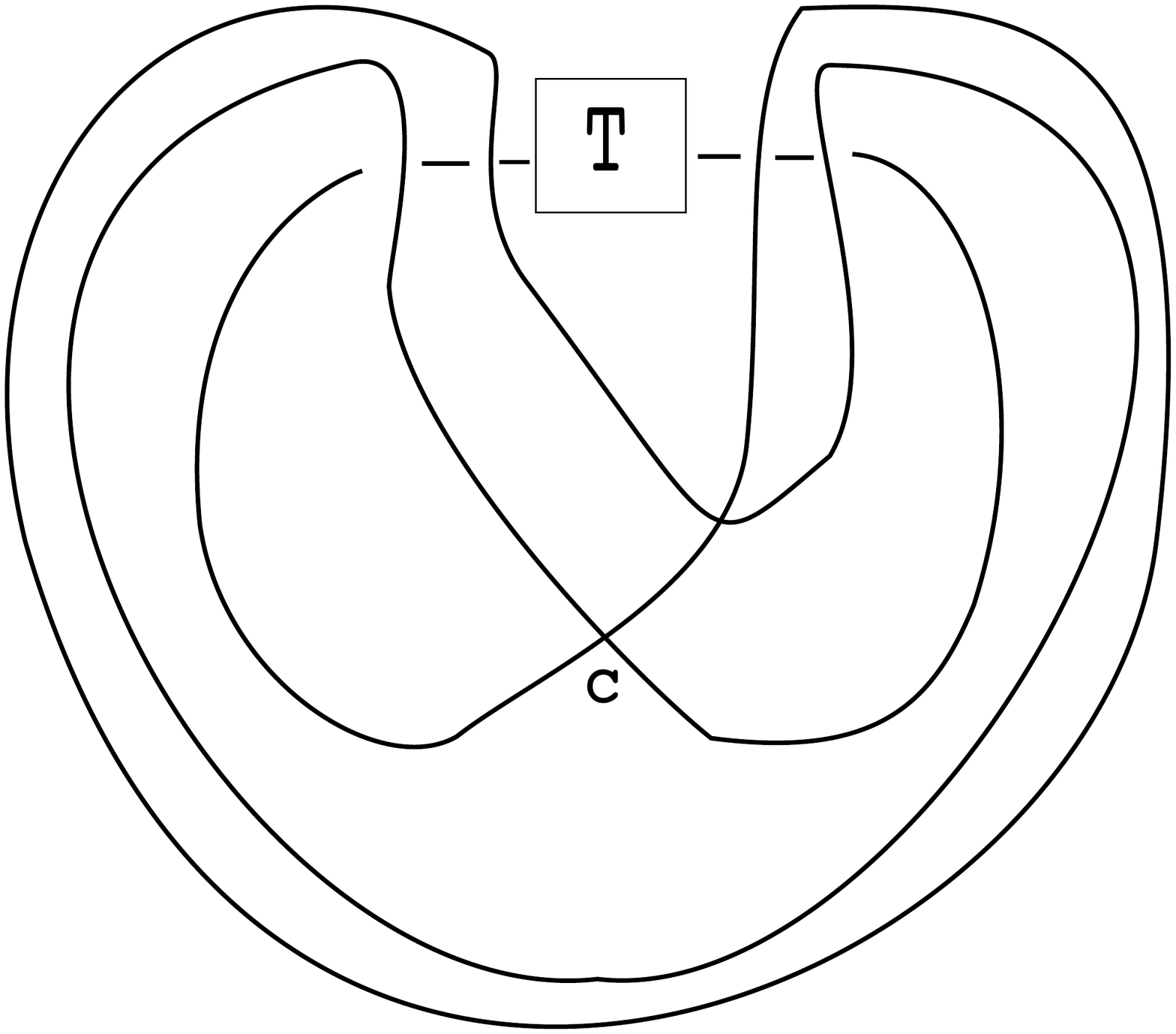}}
\caption{\small{Banana peeling:(a) shows a diagram before peeling.  The tangle {\bf T} contains the unique virtual crossing $v$.  The diagram has already been manipulated so that all kinks are concentric.
Figures (b) and (c) show the diagram after virtual detour moves on the outermost kink strand.  In the resulting diagram (c), the virtual crossing lies exterior to the strand $B$ and $B$ lies completely on top of the other strand.}}
\label{banana}
\end{figure}

Suppose instead that $v$ does not lie exterior to $B$. We will reduce this case to the previous one.    The chords not in $\mathscr{T}$ and with both endpoints on $b$ are parallel.  It follows that we can resolve crossings in $\mathscr{T}$ along $B$ so that, by a a sequence of classical Reidemeister moves, $B$ is a collection of concentric kinks as shown in (a) of Figure \ref{banana}.

Starting with the outermost kink, we consecutively perform a virtual detour move on each kink  that lies exterior to $v$, as illustrated in (b) and (c) of Figure \ref{banana}.  In the resulting diagram (depicted in (c) of Figure \ref{banana}, $v$ is exterior to $B$, and crossings between $A$ and $B$ are unchanged.  By the previous case, the resolution of $c$ does not affect the triviality of $S$.
\end{proof}

As an extension of this result, we propose the following conjecture.

\begin{conjecture}
For every virtual shadow $S$, $\tr{S}$ is even.
\end{conjecture}

We now turn our attention from unknotting numbers of knots and properties of the trivializing number to the topological notion of genus of a classical knot. In addition to bounding the unknotting number of a knot, the trivializing number provides an upper bound on genus.

\begin{lemma}\label{seifert}
A (classical) $n$-component link has at least $n$ Seifert circles.
\end{lemma}
\begin{proof}
This may be shown by induction on $n$. However, we will find it more enlightening to give the following braid theoretic proof.

In \cite{yamada}, Yamada showed that, for any link $L$, the braid index $\beta(L)$ is equal to the minimum number of Seifert circles in any diagram of $L$. Since a braid form representation of an $n$ component link must have at least one strand for each component, $\beta(L) \geq n$, and the lemma follows.
\end{proof}

\begin{theorem}\label{genus}
For any classical knot $K$, let $D_K$ be any diagram of $K$ and let $S_K$ be the shadow of $D_K$. Then $$\frac{\tr{S_K}}{2} \geq g(K),$$ where $g(K)$ is the genus of $K$.
\end{theorem}
\begin{proof}
It suffices to show that $$\frac{\tr{D}}{2} \geq g(D),$$ for every diagram $D$ of $K$, since $g(K) \le g_c(K) \leq \min_D g(D)$ where $g_c(K)$ is the canonical genus of $K$.

For every diagram $D$, it is well-known (see, for example, \cite{knotbook}) that $g(D)=\frac{c-s+1}{2}$, where $c$ is the number of crossings in $D$ and $s$ is the number of Seifert circles.

Therefore, to prove that $\frac{\tr{D}}{2} \geq g(D)$, it suffices to show that $ c - \tr{D} +1\leq s$. By Lemma \ref{triv=crosses}, we know that $\tr{D}$ is precisely the minimum number of chords that must be deleted from the chord diagram of $D$ in order to leave only parallel chords. Hence, $c- \tr{D}$ is the cardinality of the maximum set of parallel chords in the chord diagram of $D$. These $c - \tr{D}$ chords divide the chord diagram into $c -\tr{D} +1$ planar regions.

Smoothing at these $c -\tr{D}$ crossings then produces a $c -\tr{D} +1$ component link. Then, by Lemma \ref{seifert}, $D$ has at least $c - \tr{D} +1$ Seifert circles, as desired.
\end{proof}

In the future, we hope to resolve some of the open questions proposed in this paper. We also plan to extend this theory to include links and spatial graphs, both classical and virtual. Furthermore, it would be interesting to continue to relate concepts from pseudodiagram theory to more established knot invariants, as was begun in this final section.

\section*{Acknowledgements}
The authors took part in the SMALL Summer 2009 REU at Williams College, supported by Williams College, Oberlin College, and NSF grant DMS-0850577.  The authors would also like to thank Robert Silversmith for editing advice and a suggestion about Theorem \ref{knottingnotsmall}, and Colin Adams for much advice and support, particularly for proposing a better proof of Lemma \ref{seifert}. We would also like to thank Vassily Manturov for his valuable comments regarding Lemma~\ref{+nontriv} and a suggested proof for Lemma~\ref{ManturovProof}.



\bibliographystyle{abbrv}
\bibliography{pseudodiagrams}

\end{document}